\documentclass{article}

\usepackage{arxiv}

\usepackage[utf8]{inputenc} % allow utf-8 input
\usepackage[T1]{fontenc}    % use 8-bit T1 fonts
\usepackage{hyperref}       % hyperlinks
\usepackage{url}            % simple URL typesetting
\usepackage{booktabs}       % professional-quality tables
\usepackage{amsfonts}       % blackboard math symbols
\usepackage{amssymb}
\usepackage{nicefrac}       % compact symbols for 1/2, etc.
\usepackage{microtype}      % microtypography
\usepackage{lipsum}		% Can be removed after putting your text content
\usepackage{graphicx}
\usepackage{natbib}
\usepackage{doi}
\usepackage{amsmath, bm}
\usepackage{ellipsis}

\title{Exploration-oriented sampling strategies for global surrogate modeling: A
comparison between one-stage and adaptive methods}

\author{ \href{https://orcid.org/0000-0001-9722-6185}{\includegraphics[scale=0.06]{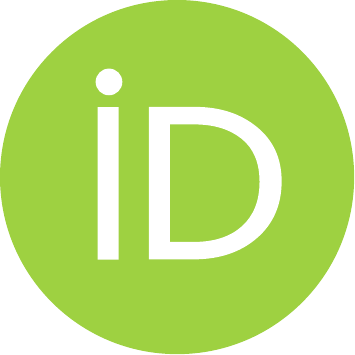}\hspace{1mm}Pietro ~Lualdi}
    \thanks{Corresponding author} \\
	Institute of Vehicle Concepts (FK)\\
	German Aerospace Center (DLR)\\
	Pfaffenwaldring 38-40, 70569, Stuttgart\\
	\texttt{pietro.lualdi@dlr.de} \\
	\And
	\href{}{\includegraphics[scale=0.06]{orcid.pdf}\hspace{1mm}Ralf ~Sturm} \\
	Institute of Vehicle Concepts (FK)\\
	German Aerospace Center (DLR)\\
	Pfaffenwaldring 38-40, 70569, Stuttgart\\
	\texttt{ralf.sturm@dlr.de} \\
	\And
	\href{}{\includegraphics[scale=0.06]{orcid.pdf}\hspace{1mm}Tjark ~Siefkes} \\
	Institute of Vehicle Concepts (FK)\\
	German Aerospace Center (DLR)\\
	Pfaffenwaldring 38-40, 70569, Stuttgart\\
	\texttt{tjark.siefkes@dlr.de} \\
}

\renewcommand{\headeright}{\textit{Lualdi et al. 2022}}
\renewcommand{\undertitle}{Preprint}
\renewcommand{\shorttitle}{\textit{}}

\hypersetup{
pdftitle={A template for the arxiv style},
pdfsubject={q-bio.NC, q-bio.QM},
pdfauthor={Pietro .~Lualdi, Ralf ~Sturm, Tjark ~Siefkes},
pdfkeywords={Exploration sampling strategies, Adaptive sampling method, Global surrogate modelling, Crashworthiness optimization},
}

\begin{document}
\maketitle

\begin{abstract}
	{Studying complex phenomena in detail by performing real experiments is often an unfeasible task. Virtual experiments using simulations are usually used to support the development process. However, numerical simulations are limited by their computational cost. Metamodeling techniques are commonly used to mimic the behavior of unknown solver functions, especially for expensive black box optimizations. If a good correlation between the surrogate model and the black box function is obtained, expensive numerical simulations can be significantly reduced. The sampling strategy, which selects a subset of samples that can adequately predict the behavior of expensive black box functions, plays an important role in the fidelity of the surrogate model. Achieving the desired metamodel accuracy with as few solver calls as possible is the main goal of global surrogate modeling. In this paper, exploration-oriented adaptive sampling strategies are compared with commonly used one-stage sampling approaches, such as Latin Hypercube Design (LHD). The difference in the quality of approximation is tested on benchmark functions from 2 up to 30 variables. Two novel sampling algorithms to get fine-grained quasi-LHDs will be proposed and an improvement to a well-known, pre-existing sequential input algorithm will be discussed. Finally, these methods are applied to a crash box design to investigate the performance when approximating highly non-linear crashworthiness problems. It is found that adaptive sampling approaches outperform one-stage methods both in terms of mathematical properties and in terms of metamodel accuracy in the majority of the tests. A proper stopping algorithm should also be employed with adaptive methods to avoid oversampling.}
\end{abstract}

% keywords can be removed

\section{Introduction}\label{sec1}
{The application of numerical simulations is prevalent in many engineering disciplines. Even if simulations cannot completely replace real-life experiments, they can considerably reduce development costs. Thus, simulations are used by engineers to gain a better understanding of problems and to identify designs with improved performance.  

In the field of Design Optimization, the best possible values of the design variables under consideration (also called “features”) are identified so that the investigated objective function(s) can be either maximized or minimized, while still satisfying all specified constraints typical of engineering problems \citep{Baier1994OptimierungID}. This can be achieved by an iterative approach using global optimization algorithms.  This approach requires a high number of evaluations of the objective and constraint functions. However, complex engineering systems such as Finite Element Analysis (FEA) and Computational Fluid Dynamics (CFD) can be computationally expensive, which makes an industrial application difficult. For instance, in \citet{Koch2008} it is reported that full car crash models in current projects at Porsche consist of about 15 million elements, leading to a simulation time of about 32 hours, using 256 processors of the latest generation. These computational costs clearly limit the application of classic optimization approaches.

In order to overcome this issue, surrogate modeling techniques are often used. Metamodels, often known as surrogate models, are mathematical approximations that can be used efficiently to deal with complex and computational expensive black box functions. Surrogate models are capable of efficiently capturing the behavior of linear and non-linear problems, by mapping the output responses of a system to the input variables.

A distinction must be drawn between global and local surrogate modelling. Local models are used to guide the optimization algorithm towards a global or local optimum. By contrast, the aim of global surrogate modeling is the creation of a model that mimics the behavior of the black box system on the entire domain. Consequently, the surrogate model can be used as an approximated replacement for the original function. This paper focuses on global surrogate modelling, and therefore the issue of local surrogate modelling is not further addressed.

Assuming that no information about the problem is known a priori (i.e., there is no knowledge about the function that correlates the input parameters with the output responses), initial designs must be selected within the design space in an effective way, to allow an efficient start of the optimization process. These initial designs are called “samples”, and together they build up the Design of Experiment (DoE). The main task of the DoE is to obtain the maximum amount of information about the unknown function under investigation. Since solver calls can be expensive, it is important to find an optimal sample distribution in the design space, in order to understand the global behavior of the model with the smallest possible number of solver calls.

This paper is structured as follows. In Sect. \ref{sec2} state-of-the-art input-oriented sampling methods are presented. Three mathematical criteria to evaluate the performance of the algorithms are introduced and sampling methods are then classified according to their adaptive nature. The potential benefits and drawbacks of such methods are described, with a particular focus on expensive black box functions such as crash simulations. An empirical approach to tune the initial parameter of a well-known adaptive strategy is then shown. Furthermore, two novel adaptive sampling strategies inspired by existing methods are presented. In Sect. \ref{sec3}, the most promising methods are tested both on a pool of benchmark functions and on a specific real crash application. Two regression models are investigated for these tests: Kriging and Support Vector Regression (SVR). Finally, Sect. \ref{sec4} completes the paper with conclusions and an outlook for future research.
$ \mathbb{R} $
}

\setlength{\ellipsisgap}{0.01em}
\section{Design of Experiment}\label{sec2}
{A DoE is a structured method for determining the relationships between input factors (independent variables) and one or more output responses (dependent variables), through the application of mathematical models. In the DoE, the input factors are systematically varied to determine their effects on the output responses, which allows the determination of the most important input factors: the identification of input factors with optimized output responses, and the interactions between input factors \citep{Fukuda2018DesignOE}. Formally, a DoE can be defined as a set of n combinations of d factor values. These combinations are usually bounded by upper and lower boundaries, so that for each independent variable $x^{k}$ it holds that $ a^k \leq x^{k}_i \leq b^k$ with $a^k, b^k \in \mathbb{R}, k = 1, 2, \ldots, d $ and $i= 1, 2, \ldots, n$ \citep{Kleijnen2004}. With no loss of generality, in this paper a design space $T = [0, 1]^d $ is considered without exception, which implies that the range of each function argument has been scaled to the unit interval and that the joint region of interest is the $k$-dimensional unit cube \citep{JOHNSON1990131}. Therefore, a DoE of size $n$ turns out to be a design composed of a set of scattered points $P = \{p_1, p_2, \ldots, p_n\}\subset[0,1]^d$ with $p_i = (p^{1}_i, p^{2}_i, \ldots, p^{d}_i),$ where the function values $Y = \{f(p_1), f(p_2), \ldots, f(p_n)\}$ are known. These scattered points $P$ are commonly known as “samples”, while the function values $Y$ are called “response values”. Finally, the data points $P$ together with their function values $Y$ are used to find the best suitable surrogate function $ \hat{f}: \mathbb{R}^d \to \mathbb{R} $ of the unknown function $f: \mathbb{R}^d \to \mathbb{R}$ which describes the mapping between inputs and outputs of the black box function. From now on, a generic $d$-dimensional dataset $P$ of $n$ sample points will be represented as a matrix of $n$ rows and $d$ columns. This matrix will be denoted by $P(n,d)$.}

\subsection{One-stage and sequential sampling methods}\label{sec2.1}
{The most common sampling strategy is the one-stage approach (also
called “one-shot” approach). These design methods consider the design
space only to generate samples, and, most importantly, to spread them
out uniformly over the entire domain. It should be noted that these
approaches do not consider the function values $Y$, since the set of datapoints $P$ is generated upfront in one stage when the response values of these samples are not yet known. It is clear that the main disadvantage of this approach is the risk of running into over- or undersampling. In oversampling, too many samples have been evaluated to achieve the desired accuracy, which can result in high computational costs. In undersampling, too few sample points may have been evaluated, which requires a new DoE to reach the expected accuracy \citep{JIANG2015532}. The reason why these methods are still widely used is their ease of implementation and their optimal coverage of the domain \citep{CROMBECQ2011683}.
An alternative sampling approach is given by sequential sampling. In
literature, sequential sampling is also known as adaptive sampling \citep{Lehmensiek2002} or active learning \citep{Sugiyama2006}. For the sake of clarity, from now on these terms will be used interchangeably. In sequential sampling, a dedicated algorithm selects a few samples (or even a single sample) and adds them
progressively into the design space until the desired accuracy is reached, or the maximum number of points allowed is exceeded. The sequential sampling approach is thus an iterative process. Unlike traditional sampling strategies, response values and samples from previous iterations can be analyzed and used to generate new samples in areas that are assumed to be the most advantageous for exploration. The ability to stop the sequential algorithm, once the desired level of accuracy has been reached, is the most significant advantage of this type of algorithm.}

\subsection{Exploration and exploitation}\label{sec2.2}
{Due to the different goals of exploration and exploitation-oriented methods, it is necessary to clarify the difference at the outset. Exploration means gaining the crucial information of the unexplored design. Peaks, valleys, sharp changes of response values, and discontinuities are just a few examples of this type of information. Since the main focus of exploration is the scanning of the design space, the response values are not considered. 
Exploitation means the investigation of a portion of the domain that has already been identified to be of interest. Refining a potential local optimum or selecting new samples in steep regions are possible exploitation-oriented sampling approaches \citep{Crombecq2009SpacefillingSD}. The main difference compared to exploration is that exploitation takes advantage of the response values of the previous iterations to select new samples. Thus, exploration and exploitation refer to input and output-based sampling, respectively. 
While it is evident that one-stage sampling approaches are exploration-oriented (also known as input-oriented), the focus of adaptive methods is rather less obvious. Adaptive methods can be exploration-oriented, exploitation-oriented, or have a hybrid focus, requiring an optimal trade-off between the goals of exploration and exploitation. 
As outlined above, the main advantage of adaptive strategies is that undersampling and oversampling phenomena can be avoided by halting the sampling process as soon as the desired accuracy is reached. Although hybrid strategies have greater potential (since they have access to more information), exploration-oriented strategies are certainly no less useful. In fact, they are particularly useful in multi-response systems where a metamodel is needed to model each response. This is a typical scenario of multi-objective or multi-constraint problems (e.g. crashworthiness optimization) where various objective and constraint functions are built on a single DoE.
The main focus of this paper is to analyze the benefits of adaptive, exploration-based strategies compared to state-of-the-art one-stage approaches.}

\subsection{Properties of sampling strategies}\label{sec2.3}
{The quality of the sample distribution over the design space can be assessed by several mathematical criteria. As proposed by \citep{JOHNSON1990131} and \citep{Morris1995ExploratoryDF} for computer experiments, at least two criteria, the space-filling, and the non-collapsing criteria must be satisfied. Another ideal property described by \citep{CROMBECQ2011683} is the granularity of the design. Granularity is often the main difference between sequential and one-stage approaches. In this section, these three properties are introduced and discussed.}

\subsubsection{Space-filling criterion}\label{sec2.3.1}
{The space-filling criterion describes the uniformity of the sampling distribution in the domain. Since no details about the functional behavior of the design parameters are available, it is crucial to gain information from the entire design space. In that sense, the design has to be “space-filling”, which implies that its samples have to be evenly spread over the whole domain. To describe mathematically the even distribution of a design set $P$, it is necessary to define a mathematical criterion to be optimized. The most widely used are the $L_2$ norm (or maxmin) and the  $\phi_p$ criterion:}

\begin{equation}
\label{eq1}
\min_{p_i,p_j \in P} \quad \sqrt{\sum_{k=1}^{d}{(p^k_i-p^k_j)^2}}
\end{equation}

\begin{equation}
\label{eq2}
 (\sum_{p_i,p_j \in P}^{}{[\sum_{k=1}^{d}{(p^k_i-p^k_j)^2}]^{-p}})^{1/p}
\end{equation}

{
Regardless of the choice of criterion, the higher the value, the better the space-filling properties. Because of numerical instability issues (especially for large values of the power $p$) of the $\phi_p$ criterion, the maxmin criterion ($L_2$  norm) is used in this work to compare different designs based on their space-filling properties.

In order to measure the relative isolation of a point $p$ from the existing dataset $P$, the Crowding Distance Metric (CDM) can be used. This distance metric deployed first by \citep{zhang2012adaptive} and resumed by \citep{GARUD2017103} can be used to place new points in relatively unexplored regions and as far away from the existing points as possible. Its mathematical formulation, based on Euclidean norm, is given in Eq. \ref{eq3}:
}

\begin{equation}
\label{eq3}
CDM(P, \bm{p}) = \sum_{i=1}^{n}{(\|\bm{p}-\bm{p_i}\|)^2}
\end{equation}

{Again here, the greater the CDM measure, the better the isolation.}

\subsubsection{Non-collapsing criterion}\label{sec2.3.2}
{For an optimal design, non-collapsing properties are beneficial too. This feature, also called the projective property, is guaranteed when each coordinate $p^k_i$ of every sample $\mathbf{p_i} \in P$ is strictly unique. To understand this concept, it is useful to consider a problem where one input variable has almost no influence on the response function values. If two samples differ only in this variable, they will “collapse”, i.e., they can be considered as the same point evaluated twice \citep{Husslage2006MaximinDF}. In the context of expensive black box simulations, this is not desirable since it leads to unnecessary computational costs. Mathematically, to describe the quality of a design in terms of non-collapsing property, the minus infinity norm is used to define the \it{minimum projected distance}:}

\begin{equation}
\label{eq4}
\min_{p_i,p_j \in P}\|\bm{p_i}-\bm{p_j}\|_{-\infty}=\min_{p_i,p_j \in P}\min_{1\leq k \leq d}|p^k_i-p^k_j|
\end{equation}

{Similar to the maxmin criterion, the higher the minimum projected distance, the better the non-collapsing property of a specific design.}

\subsubsection{Granularity}\label{sec2.3.3}
{The granularity of the strategy plays a major role in the quality of the sampling design. Two main aspects of granularity should be highlighted. To begin with, an optimal granularity belongs to a sampling strategy for which it is not necessary to know the total number of samples in advance. Hence, classic one-stage sampling approaches have the worst possible granularity. Moreover, an ideal design should be fine-grained. This means that, considering a sequential sampling strategy, the best possible scenario is if one point is selected and added to the design at the end of each iteration. By contrast, a coarse-grained sequential sampling strategy selects large batches of design points after every iteration, which can lead to under- and oversampling.}

\subsection{Overview of existing methods}\label{sec2.4}
\begin{figure}[!htb]
\minipage{1.0\textwidth}
\centering
  \includegraphics[]{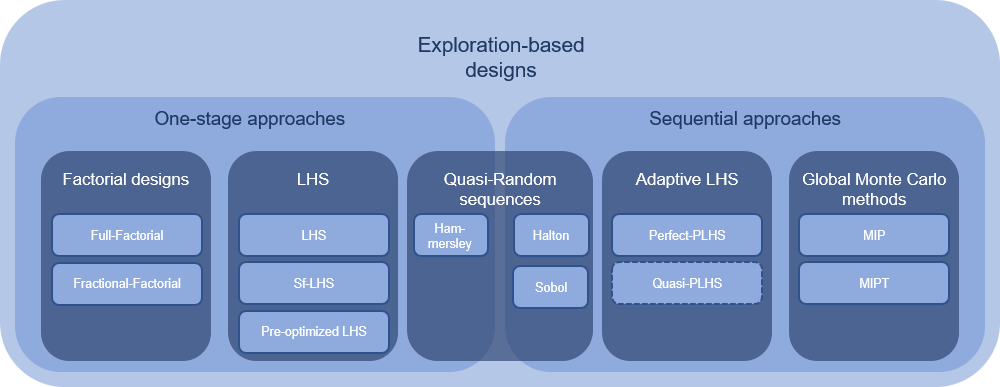}
  \caption{Classification of exploration-based sampling methods}
  \label{fig_class}
\endminipage\hfill
\end{figure}

{In Figure ~\ref{fig_class}, common DoE types and some recently published input-based sampling algorithms are classified according to their adaptive nature (one-shot or sequential approach). The goal is to distinguish methods that can be employed as sequential sampling strategies, i.e. to add new points to the initial dataset from “static” one-shot approaches. After a brief description, the most promising methods are tested and compared for different non-linear optimization problems.} 

% test picture

\subsubsection{Factorial designs}\label{sec2.4.1}
{
Factorial designs represent a category of space-filling-oriented sampling strategies. Full-factorial designs are $d$-dimensional grids of points evenly distributed over the domain. Having d factors (i.e., the dimension of the problem) and an integer $m$, which defines the number of levels the $d-th$ level is divided into (i.e., the resolution of each input variable), a full-factorial design is given by $m^d$ samples, where each sample is a specific combination of the $d$ factors \citep{Draper1996}. This type of DoE results in the best space-filling performance achievable, owing to the even distribution of samples in the design space. Under the assumption that the input variables are equally important, a space-filling design would be an ideal sampling strategy in terms of design quality. However, full-factorial designs are extremely expensive in terms of the function evaluations required. The number of samples grows exponentially with the number of dimensions, which can be unacceptable in the field of FEM and CFD simulations. Additionally, a full-factorial strategy turns out to be a coarse-grained strategy. The total number of samples must be known in advance, which makes this strategy unfeasible for the sampling refinement typical of adaptive strategies. Furthermore, such designs are prone to the danger of aliasing, as discussed by \citep{Wu2017}. Moreover, A full factorial design has the worst possible non-collapsing properties. By definition, the samples are not strictly unique, since they share their coordinates over the $d$-dimensional grid. A partial solution to handle this drawback is given by fractional designs. Starting from a full-factorial grid, the fractional designs solution brings improvements in terms of total runtime and projective properties by removing points. Nevertheless, given their previously highlighted disadvantages, factorial designs will not be considered in this study.
}

\subsubsection{Latin Hypercube Sampling}\label{sec2.4.2}
{
Latin Hypercube Sampling (LHS), also called Latin Hypercube Designs (LHD), are extremely useful designs of experiment in the field of black box optimization \citep{vanDam2007}. Their well-understood mathematical properties, ease of implementation, and speed make them widely used sampling methods \citep{CROMBECQ2011683}.
For a unit cube $T$ in a $d$-dimensional space $T=[0,1]^d$ divided into $n$ intervals (where $n$ is the sample size) with an equal length of $1/n$   along each axis, LHS creates $n$ equally probable intervals indexed by $q = 1,\ldots,n$ corresponding to $[0,1/n),[1/n,2/n),\ldots,[(n-1)/n,1]$ for each dimension. LHS can be represented as an $n$-by-$d$ sample matrix $[x_{i,j}]$ $(i = 1,\ldots ,n; j = 1,\ldots,d)$, where $x_{i,j} \in [0,1]$ such that $x_{i,j}$ in the $j-th$ column belongs to only one interval. Therefore, $q$ is a random permutation of ${1,2,\ldots,n}$ for each column and each row of the matrix. Given a new set of binary variables $y_{q,j}$ such that
}

\begin{equation}
\label{eq5}
    y_{q,j}=\left\{
                \begin{array}{ll}
                  1 \;\;\;\; \textnormal{ if there exist any $i$ for which $x_{i,j}$  lies in the interval $q$}\\\
                  0 \;\;\;\; \textnormal{otherwise}\
                \end{array}
              \right.
\end{equation}

{
From this follows that a dataset $P(n,d)$ is called a Latin Hypercube (i.e. $LHS(n,d)$) if the following condition is met:
}

\begin{equation}
    \label{eq6}
    \frac{\sum_{j=1}^{d}{\sum_{q=1}^{n}{y_{q,j}}}}{n \cdot d}=1
\end{equation}

{
The left-hand side of equation \ref{eq6} yields a measure of how close a dataset is to an LHS. This value can range between 1/n (this happens when all the samples lie in the same interval at every dimension) and 1 (when the dataset is an LHD).
}

{
According to the given mathematical definition and assuming a uniform distribution of the samples along every dimension, a sample is only “Latin Hypercube” if it possesses one-dimensional projection properties. Such a sample, however, is only guaranteed to maximize its non-collapsing properties, while the space-filling properties are not necessarily accounted for \citep{SHEIKHOLESLAMI2017109}. 
LHS with better space-filling properties can be obtained with maxmin LHS, which means a Latin Hypercube Design where the minimum distance between two samples has been maximized. Several distance measures such as $l^\infty,l^1,l^2$ could be considered in the maxmin criterion. This paper focuses on the Euclidean distance since its definition is consistent with the formula of the space-filling criterion. A valuable database of pre-optimized LHSs has been published by \citet{GROSSO2009541}, \citet{Husslage2006MaximinDF}, and \citet{vanDam2007} and it is available at \url{https://spacefillingdesigns.nl}. It has to be noted that these final designs achieved an optimal or semi-optimal solution by placing the samples in the middle of their intervals. In this way, the projected distance results are uniform along each direction. However, given $n$ and $d$, $(n!)^{d-1}$ different LHDs can be generated. Therefore, since these LHDs cannot be generated in real-time (each design requires hours of optimization) this database is directly used within this study. Since the number of samples has to be known a priori, these designs are not applicable for adaptive approaches.
}

\subsubsection{Quasi-Random sequences}\label{sec2.4.3}
\begin{figure}[!htb]
\minipage{1.0\textwidth}
\centering
  \includegraphics[width=0.49\columnwidth]{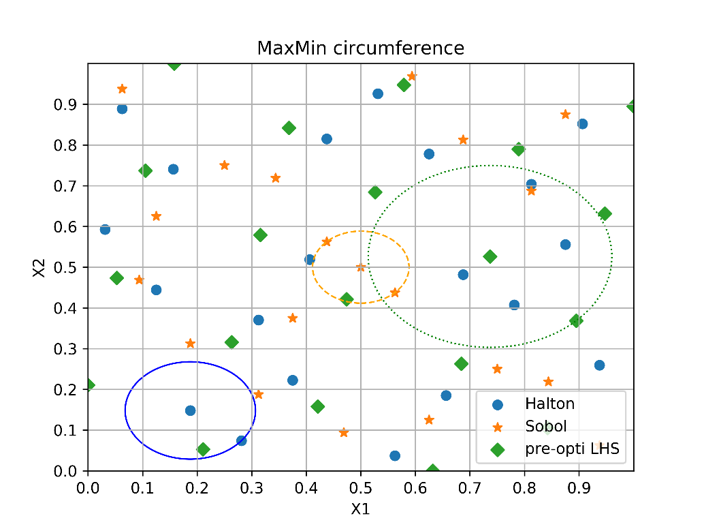}
  \caption{Comparison of low-discrepancy sequence (Halton and Sobol) with a pre-optimized LHS (20 samples). The minimum distance for each design is highlighted as a circumference with radius equal to the Euclidean distance between the two closest samples}
  \label{fig_02}
\endminipage\hfill
\end{figure}

{
Quasi-random sequences, also called “low-discrepancy sequences”, are placed between one-stage and sequential methods in the proposed classification. Due to their adaptive nature, some correspond to sequential approaches (e.g., \citep{SOBOL196786}, \citep{Halton1964Algorithm2R}), while others do not (e.g. \citep{Hammersley1964}). All of them, however, are deterministic strategies, which means that they utilize deterministic routines designed for space-filling goals. To achieve such mathematical properties, a low discrepancy criterion is adopted, depending on the strategy. The discrepancy metric, deﬁned by Ilya M. Sobol, is the maximum deviation between the theoretical density $dt=1/n$   and the point density $d_i$ in an arbitrary hyper-parallelepiped ($T_i$) within the parameter space (hypercube) \citep{SALTELLI2010259},\citep{Burhenne2011SAMPLINGBO}. A low discrepancy of the design is therefore guaranteed when the two densities are close to proportional.
}
{
The projective properties of quasi-random sequences are the subject of debate; some authors assess them as average quality \citep{CROMBECQ2011683}, while others say that they are generally poor \citep{SHEIKHOLESLAMI2017109}. One point on which these authors agree is that unwanted correlation between the input variables might arise, especially in high dimensional spaces. Furthermore, their space-filling properties are poor for a relatively small number of samples. This can be easily shown by comparing those methods using maxmin radius of a (hyper) sphere with a space-filling pre-optimized LHS design. The  maxmin radii for the three methods for a 2-dimensional space with 20 samples are depicted in Figure \ref{fig_02}.
}

{
A further helpful visualization tool (which can also be used for space-filling sampling algorithms) is the description of sampling density using the Voronoi tessellation. Given a set of points $X \subset {\mathbb{R}}^d$, for any point $x_i \in X$ (denoted as "generators”), the Voronoi region (or Voronoi cell) $\hat{V}_i \in {\mathbb{R}}^d$ contains all the points belonging to the domain which lie closer to $x_i$ than to any other generator in $X$. The Voronoi tessellation is given by the complete set of Voronoi regions $\{\hat{V}_1,\hat{V}_2,\ldots,\hat{V}_n\}$ which tessellate the whole domain. Formally, a Voronoi and a generic tessellation are defined as follows. Given a set of points $ \{x_i \}_{i=1}^n=\{x_1,x_2,\ldots,x_n \} $ belonging to the closed set $\bar{\Omega} \in {\mathbb{R}}^n$, the Voronoi region $\hat{V}_i$ corresponding to the point $x_i$ is defined by:
}

\begin{equation}
    \label{eq7}
    \hat{V}_i= \{x \in \Omega \mid |x-x_i|<|x-x_j | \;\; for \;\; j=1,\ldots,n,\;\; j \neq i \}
\end{equation}

{
Given an open set  $ \Omega \in \mathbb{R}^d $ the set ${V_1,V_2,\ldots,V_n}$ is called a tessellation of $\Omega$ if $V_i \subset \bar{\Omega}$ for $i=1,\ldots,n, V_i \cap V_j = \emptyset$ for $i \neq j$, and $ U_{i=1}^n V_i=\Omega$ \citep{Du1999}, \citep{Burns2009CENTROIDALVT}.
}

\begin{figure}[!htb]
\minipage{1.0\textwidth}
\centering
  \includegraphics[]{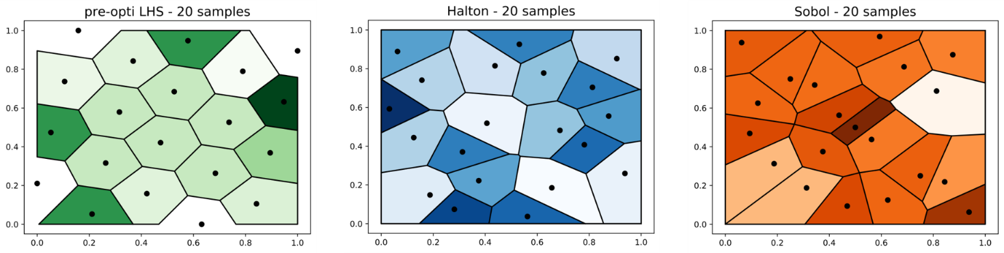}
  \caption{Voronoi tessellation for a pre-optimized LHS, a Halton, and a Sobol design respectively.}
  \label{fig_03}
\endminipage\hfill
\end{figure}

{
In Figure \ref{fig_03} the Voronoi tessellation is shown for the same sampling distribution depicted in Figure \ref{fig_02} to emphasize its subpar space-filling properties. The darkness of the cells refers to their area. The darker the cell, the less optimal its sampling position. The colors in Figure \ref{fig_03} are normalized between the largest and the smallest Voronoi region in each diagram. The diagrams show that Sobol and Halton designs do not guarantee an ideal uniform space distribution such as a pre-optimized LHS.
}
{
Nevertheless, Sobol and Halton sequences, due to their adaptive nature and fine-grained refinement, will be considered in this study. Their implementation is already available in the numerical tool ChaosPy \citep{FEINBERG201546}, which is used to assess the performance of these methods.
}

\subsubsection{Adaptive Latin Hypercube Sampling}\label{sec2.4.4}
{
To make the approach suitable for a sequential strategy, several attempts to modify the original one-stage algorithm of the LHS have been made. Husslage presented a method to sequentially generate LHDs \citep{Husslage2006MaximinDF}. Crombecq modified it to a less coarse-grained version: Nested Latin Hypercubes \citep{CROMBECQ2011683}. Sheikholeslami also presented the Perfect Progressive Latin Hypercube Sampling (perfect-PLHS) and the Quasi Progressive Latin Hypercube Sampling (quasi-PLHS) \citep{SHEIKHOLESLAMI2017109}. The main idea behind these methods is the sequential LHS generation (i.e., without re-building a new design) by adding new slices (or layers) to an already generated LHS. To ensure that the sampling algorithm remains an LHS, a refinement of the starting grid is done at each iteration. According to the chosen method, at the $i-th$ iteration $k_i$ new points ($k_i=n_{i-1}-1$  for Nested Latin Hypercubes and $k_i=2 \cdot n_{i-1}$ for perfect-PLHS) are selected and added to the starting design set. Despite the above-average projective properties, there are two main issues related to sequential LHS. The first one is related to granularity: one-grained and fine-grained strategies are inoperable to ensure the designs remain complete LHSs. As explained above, this is not ideal for time-consuming simulations, since it can lead to oversampling. Secondly, sequential LHSs are more likely to get stuck in local optima. Even if a point is optimally chosen in terms of space-filling and non-collapsing properties at the $i-th$ iteration, it is not guaranteed that the applied sampling can provide better designs (in terms of the mathematically defined criteria) at the next iterative step. In this case, the sampling is trapped in a local optimum. Due to this limitation, and the expected high number of required sampling points, the mentioned sequential LHS are not suitable for expensive black box functions and will not be further investigated in this study.
}

\subsubsection{Global Monte Carlo methods}\label{sec2.4.5}
{
Monte Carlo algorithms use the process of repeated random simulations to estimate unknown parameters necessary to improve response values. Regarding adaptive sampling, these methods are often used to discretize the domain (or portions of it) by means of many candidate points. Thus, given a design space containing infinite points, the problem is simplified by considering a finite set. These points are then compared based on a given criterion, and a selection is used for defining the next sample location(s). These random points might be very useful for making predictions or gathering information about the design space. Below, two state-of-the-art sequential Monte Carlo methods published by \citep{CROMBECQ2011683} are investigated.
\\
The first algorithm is the mc-inter-proj (MIP), which aims at maximizing the objective function in Eq. \ref{eq8}:
}

\begin{equation}
    \label{eq8}
    MIP(P,\bm{p}) = \frac{\sqrt[d]{n+1}-1}{2} \min_{p_i\in P}\|\bm{p_i}-\bm{p}\|_{2} + \frac{n+1}{2} \min_{p_i\in P}\|\bm{p_i}-\bm{p}\|_{- \infty}
\end{equation}

{
where $P$ is the set of the previously evaluated samples, $\mathbf{p}$ is a new candidate point, $d$ is the dimension of the problem, and $n$ is the size of the design, respectively. At each iteration, depending on the current number of samples $n$,$100n$ random candidate points are generated. Among them, the point which maximizes the objective function is selected as the new sample to be included in the sample set $P$. Note that both the definition of the space-filling criterion (Eq. \ref{eq1}) and the definition of the non-collapsing criterion (Eq. \ref{eq4}) affect the outcome of this formula and are scaled according to $n$ and $d$.
\\
\\
A more efficient and improved version of this function is obtained by replacing the projected distance function with a threshold function. Since such an objective function is challenging to maximize, the main idea behind the mc-inter-proj-th (MIPT) method is the simplification of this complex function by discarding candidates which lie too close to each other. The remaining points are then ranked on their intersite distance. Therefore, the threshold version of this Monte Carlo method changes as follows:
}

\begin{equation}
    \label{eq9}
    MIP(P,\bm{p}) = \left\{
                        \begin{array}{ll}
                         0 \;\;\;\;\;\;\;\;\;\;\;\;\;\;\;\;\;\;\;\;\;\;\;\;\;\;\;\;\;\;\;\; \textnormal{ if } \min_{p_i\in P}\|\bm{p_i}-\bm{p}\|_{- \infty}<d \\
                         \min_{p_i\in P}\|\bm{p_i}-\bm{p}\|_{2} \;\;\;\;\; \textnormal{ if }
                         \min_{p_i\in P}\|\bm{p_i}-\bm{p}\|_{- \infty} \geq d
                        \end{array}
                        \right.
\end{equation}

{
where the threshold $d_{min}$ is defined by a tolerance parameter $\alpha$ in which has a domain of $[0,1]$:
}
\begin{equation}
    \label{eq10}
    d_{min} = \frac{2\alpha}{n}
\end{equation}

{
The tolerance parameter $\alpha$ defines the balance between the space-filling and non-collapsing properties. Low values of $\alpha$ lead to a reduction of the projected distance constraint. Therefore, fewer candidates are discarded. On the other hand, high values of $\alpha$ result in a strict constraint to be satisfied. This reduces the chance of finding a valid candidate. 
}

\begin{figure}[!htb]
\minipage{1.0\textwidth}
\centering
  \includegraphics[]{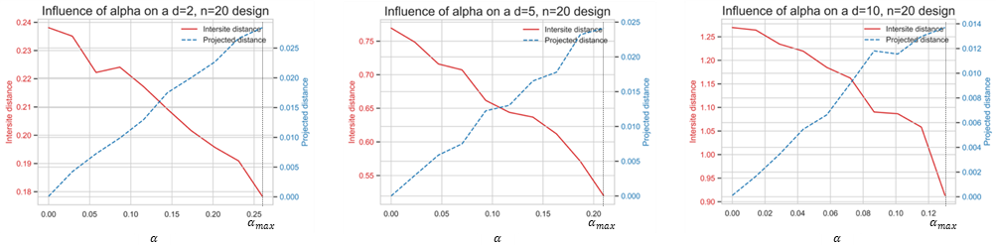}
  \caption{Influence of the tolerance parameter $\alpha$ on intersite – and projected distance for a two, five, and ten-dimensional problem with 20 samples.}
  \label{fig_04}
\endminipage\hfill
\end{figure}

{
As shown in Figure \ref{fig_04}, the choice of an optimal $\alpha$ is considerably affected by the dimensionality of the problem. Therefore, it should be adjusted to fit the problem being investigated to achieve optimal performance. A simple approach to handle this issue will be proposed in Sect. \ref{sec2.5.1}
}

\subsection{Proposed approaches}\label{sec2.5}
{
In this section, a simple empirical method to mitigate the issue related to the tuning of $\alpha$ in MIPT is presented. Scaling this parameter accordingly to the dimensionality of the problem is necessary to ensure the best trade-off between space-filling and projective properties.
\\
\\
Furthermore, two new adaptive sampling methods will be proposed in this section. As explained in the dedicated section, LHDs have well-known mathematical properties. The goal of these two proposed methods is to provide two adaptive algorithms that can resemble an LHS as much as possible, but with optimal granularity. Additionally, unlike MIPT, it is intended to provide methods that are completely independent of parameter choice. The proposed approaches are called Fluttering perfect-Progressive Latin Hypercube Sampling (FpPLHS) and Monte Carlo quasi-Latin Hypercube Sampling (MqPLHS).  These methods are inspired by the perfect-PLHS and quasi-PLHS respectively, presented by \citep{SHEIKHOLESLAMI2017109}.
It is worth pointing out from the beginning that both of these methods cannot be classified as LHS because they do not guarantee the dataset at each iterative step to be an LHS. However, they tend to approach the properties of an LHD as closely as possible and both guarantee the granularity of one.
}

\subsubsection{Tuning of parameter $\alpha$ for MIPT}\label{sec2.5.1}
{
In this section, the choice of the parameter $\alpha$ is further investigated and an empirical method for its selection is proposed. Although this issue has not been discussed in detail by its authors, the choice of the $\alpha$-parameter is directly influenced by the dimensionality of the problem. There are two main reasons for why $\alpha$ needs to be tuned. 
Firstly, according to Eq. \ref{eq9} and Eq. \ref{eq10}, a careless choice of $\alpha$ could lead to all points being rejected. Excessively high values of $\alpha$ aim to find candidates with optimal projective properties. However, keeping in mind that points are randomly generated, there is no guarantee that there is at least one point that will satisfy this condition.
\\
The second reason why proper tuning of $\alpha$  is important is given in Figure \ref{fig_04}. As discussed in the previous section, the ideal trade-off between projective and space-filling properties varies depending on the dimensionality of the problem. The intuition of why $\alpha$ depends on the dimensionality is that the potential for intersite- and projected distance increases with the dimension. The threshold intersite distance depends on $\alpha$, therefore $\alpha$ has to depend on the dimensionality to ensure that intersite- and projected distance have consistent influence on the point selected as the dimension changes.
\\
\\
The intuition behind adjusting the $\alpha$ parameter is self-evident in Figure \ref{fig_04}. Despite the differences found when varying the dimensionality, the trend of both the projective- and intersite distance looks consistent. In particular, in each of the three cases shown, the ideal trade-off between the non-collapsing and space-filling properties is roughly in the middle of the x-axis (i.e. $\alpha$). More precisely, the optimal value of $\alpha$ seems to be halfway between 0 and the maximum $\alpha$ value (hereafter called $\alpha_{max}$ for simplicity) such that at least one candidate point satisfies the threshold of Eq. \ref{eq9}. This suggests that $\alpha_{max}$  should first be calculated and then $\alpha$  should be set equal to $\alpha_{max}/2$. To calculate $\alpha_{max}$, it is necessary to find the candidate point that guarantees the maximum projective distance first. By reversing Eq. \ref{eq10} and replacing $d_{min}$ with the maximum projective distance calculated from the candidate points, $\alpha_{max}$ is obtained.
}

\begin{figure}[!htb]
\minipage{1.0\textwidth}
\centering
  \includegraphics[width=0.9\columnwidth]{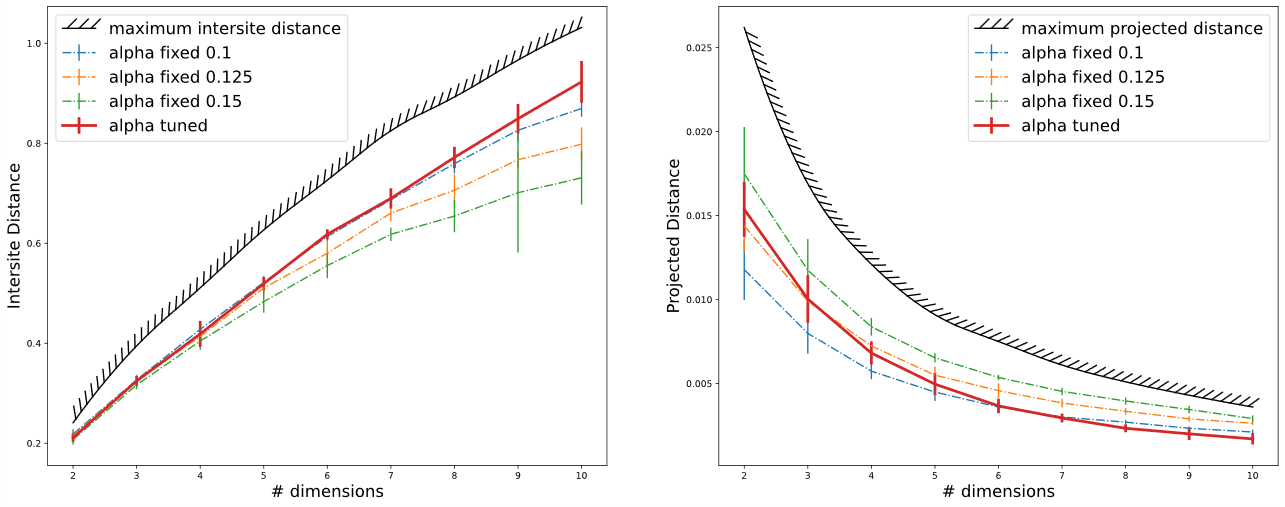}
  % \includesvg{pics/fig_05.svg}
  \caption{Comparison of the MIPT strategy with automatically adjusted alpha and fixed alpha: Intersite- (left) and projective distance (right).}
  \label{fig_05}
\endminipage\hfill
\end{figure}

{
Figure \ref{fig_05} illustrates the improvements brought by the proposed empirical approach. The conventional MIPT approach with $\alpha$ fixed (values 0.1, 0.125, and 0.15 are considered) is compared to the modified approach where $\alpha$ is automatically adjusted (“alpha tuned” in Figure \ref{fig_05}). The value of the intersite- and projective distance is shown as the dimensionality of the problem changes. In both graphs, starting from an optimal dataset of size $10 \cdot d$, a maximum estimate of the two values is represented by the path with a ticked style in black. This line was estimated using the maximum value among $1000 \cdot d$ random points in the domain. In other words, this line marks the boundary between the feasible (below the line) and unfeasible region (above the line) given by purely geometrical constraints. In both graphs, the solid red line looks like an optimal solution as dimensionality varies. This line stays relatively close to the maximum estimate, while maintaining a certain safety margin. On the contrary, observing the graph of the projective distance, it is clear that the black curve decreases faster than the three curves with $\alpha$ fixed. This will unavoidably lead these three curves to cross the unfeasible region as dimensionality increases. As for the intersite distance, better or comparable performance is achieved on average by the automatic adjustment of $\alpha$. For the above reasons, the approach proposed here seems to optimally balance projective- and space-filling properties and will therefore be used for further testing.
}

\subsubsection{Fluttering perfect-Progressive Latin Hypercube Sampling (FpPLHS)}\label{sec2.5.2}
{
The first algorithm presented here stems from the main idea of the perfect-PLHS. This approach guarantees an LHS at each iteration according to a doubling procedure. This means that, starting from an initial sample size with $P_0$ with $n_0$ samples, the first iteration will add a new LHS with $n_0$ samples to the existing dataset. Similarly, in the next iteration $2 \cdot n_0$ new samples will be added, bringing the dataset $P_2$ to have $n_2=4n_0$ total samples. This implies that the size of the dataset $P_j$ grows geometrically as $n_0 \cdot 2^j$. Due to the coarse granularity, this approach is not suitable for investigating expensive black box functions such as FEM or CFD simulation. Therefore, with the help of the example shown in Figure \ref{fig_06}, the steps of the FpPLHS approach are illustrated below.
}

\begin{description}
  \item[$ \cdot $ Step 1:] Divide the domain into equal intervals according to the number of samples $n_0=4$
  \item[$ \cdot $ Step 2:] Split each interval equally to get $2 \cdot n_0=8$ intervals along each variable
  \item[$ \cdot $ Step 3:] Prune the intervals that are already covered by existing samples along each variable
  \item[$ \cdot $ Step 4:] Generate at least $10 \cdot n_j$ independent space-filling LHS (or slices) in the residual domain gaps
  \item[$ \cdot $ Step 5:] Take the slice that, when added to the existing dataset, maximizes the space-filling properties according to the formula Eq. \ref{eq1}
  \item[$ \cdot $ Step 6:] Sort the points of the added slice according to the Crowding Distance Metric (Eq. \ref{eq3}). Each point of the slice will be added to the dataset according to this final order.
\end{description}

{
To ensure a one-grained strategy, the dataset “flutters” between one LHS and the next one, with intermediate steps not guaranteed to be an LHS. Compared to perfect-PLHS, in this approach, the good projective properties of LHS are balanced by the improvement of space-filling properties ensured by steps 4-5. Step 6 helps to avoid populating already densely populated regions and avoid undesired correlations with existing samples. Finally, by considering $10 \cdot n_j$ Latin Hypercube slices at each new iteration, a fairly efficient approach even for large datasets in high dimensional problems is ensured. Namely, once the slice that maximizes the space-filling properties is identified, sorting the points according to CDM is a computationally inexpensive operation.
}

\begin{figure}[!htb]
\minipage{1.0\textwidth}
\centering
  \includegraphics[width=0.8\columnwidth]{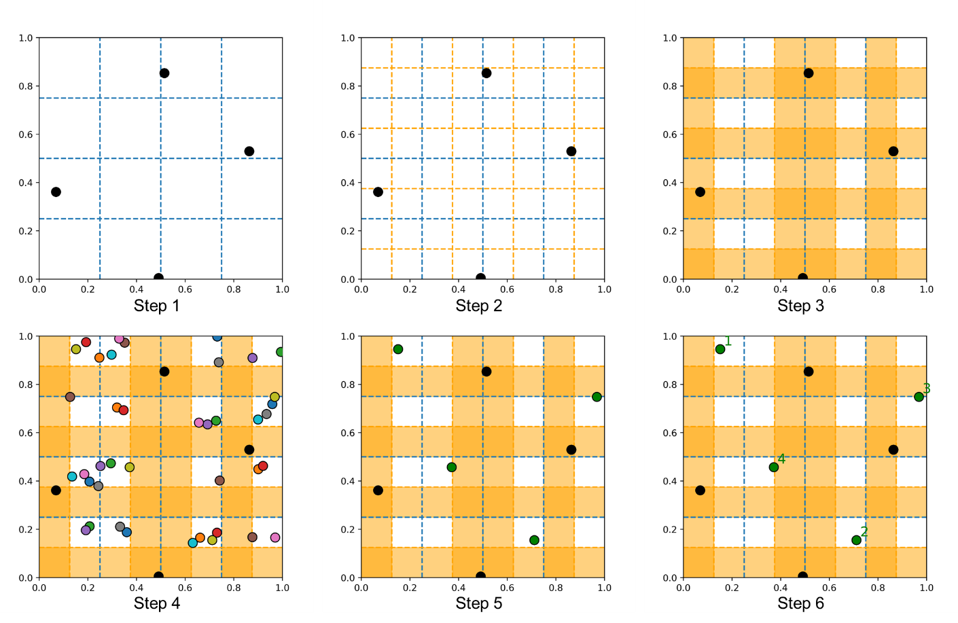}
  \caption{Iteration steps of the FpPLHS method for a dataset with $n_0=4$ and $d =2$}
  \label{fig_06}
\endminipage\hfill
\end{figure}

\subsubsection{Monte Carlo quasi-Latin Hypercube Sampling (MqPLHS)}\label{sec2.5.3}
{
A Monte Carlo approach to the quasi-PLHS is presented in the following. The main idea of this algorithm is based on the left side of Eq. \ref{eq6}. Indeed, the generation of sequential Latin Hypercube can be seen as the following optimization problem:
}

\begin{equation}
    \label{eq11}
    f_{lhs} \:\:\:\:\:\: max\Bigg( \Bigg.\frac{\sum_{j=1}^{d}{\sum_{q=1}^{n}{y_{q,j}}}}{n \cdot d} \Bigg. \Bigg)
\end{equation}

{
Starting from an initial Latin Hypercube $LHS(n,d)$, $n \cdot 100$ random points are generated. When possible, it would be recommended to start with an optimized LHS if one is available. For each of the generated candidate points, the objective function of Eq. \ref{eq11} is evaluated assuming this candidate is added to the initial dataset. Multiple candidate points can produce the same result for the objective function.  Therefore, a list containing all points that maximize the objective function is returned. Among the remaining candidates, the one that, combined with the initial dataset, guarantees the best space-filling properties according to Eq. \ref{eq1} is finally chosen. The pseudo-code of MqPLHS is described in detail below.
}

\begin{figure}[!htb]
\minipage{1.0\textwidth}
\centering
  \includegraphics[width=0.85\columnwidth]{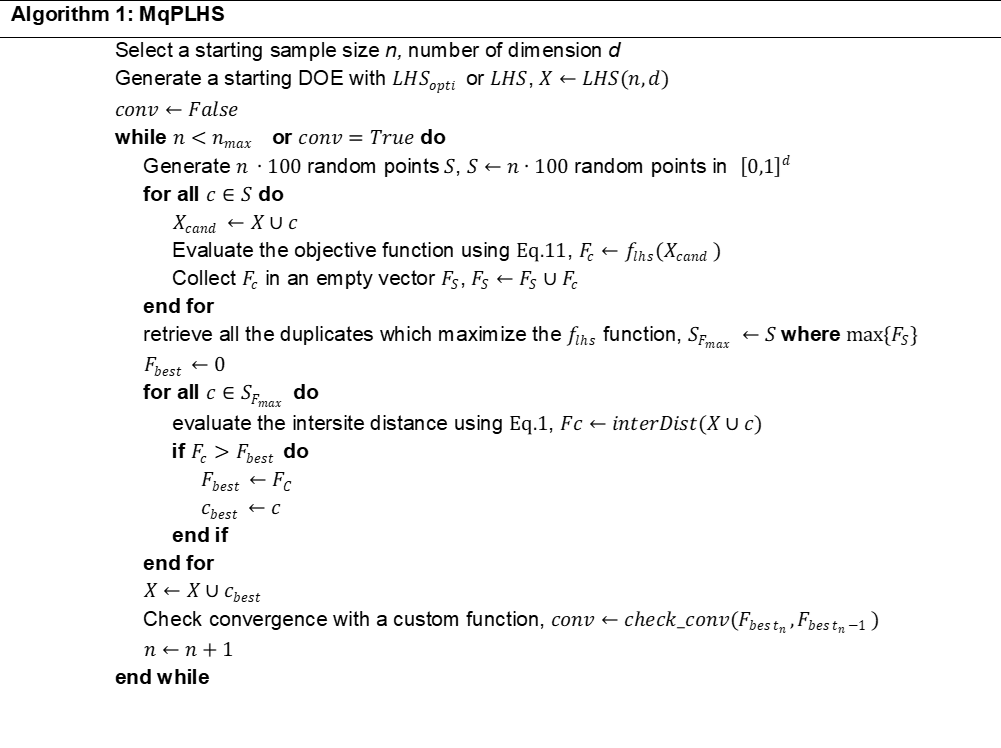}
  \label{algo_02}
\endminipage\hfill
\end{figure}

{
It is reasonable to expect that there will be several duplicates that maximize the objective function at each step. As a matter of fact, there may be more than one ideal interval to add a point to so that the design approaches the properties of a Latin Hypercube Design as closely as possible. This step is key in reducing the candidate pool and ensuring that the one that maximizes space-filling properties is chosen efficiently. 
\\
\\
It is worth mentioning that this method is essentially a greedy heuristic algorithm. This means that the choice of the ideal candidate maximizes the objective function only locally at a given iteration step. Unlike method MqPLHS, the global maximum, i.e. having the objective function equal to one, is not guaranteed to be reached.
}

\begin{figure}[!htb]
\minipage{1.0\textwidth}
\centering
  \includegraphics[width=0.85\columnwidth]{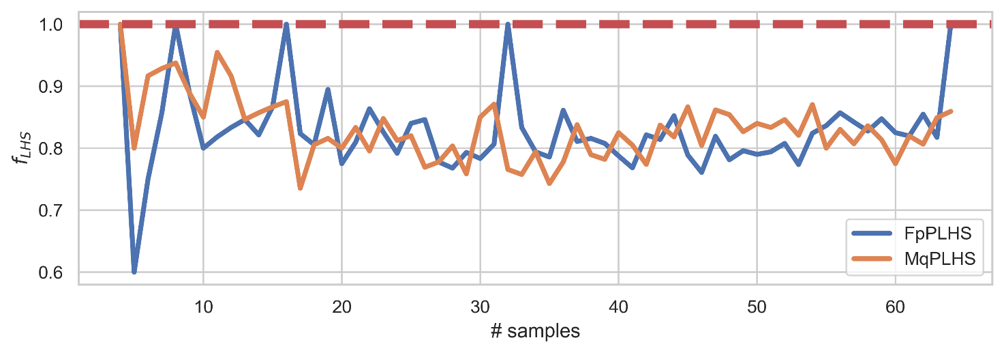}
  \caption{Comparison of the objective function $f_{lhs}$ value between FpPLHS and MqPLHS}
  \label{fig_07}
\endminipage\hfill
\end{figure}

{For verification, classic and pre-optimized LHSs are going to be investigated as one-shot sampling approaches. Concerning adaptive methods, Sobol, Halton, MIP, and MIPT will be considered.}

\section{Verification and discussion}\label{sec3}
\subsection{Evaluation of the results}\label{sec3.1}
{The sampling methods presented in the previous section are evaluated in terms of space-filling, non-collapsing, and granularity properties. For this preliminary testing phase (Sect. \ref{sec3.4}), three different design spaces are considered (2-,5- and 10-dimensional domain), to investigate the influence of the dimensionality of the problem and size of the dataset on the performance of each algorithm. The size of the DoE is set to 10 samples per input parameter ($n=10 \cdot d$) as recommended by \citep{Schreiber2006} and \citep{Jones1998EfficientGO}. The MIP, MIPT, FpPLHS and MqPLHS sampling methods start from a 10-sample, pre-optimized Latin Hypercube design. Each test is repeated 30 times to identify 95\% confidence interval bandwidths due to their stochastic nature. The most promising methods are then tested on several benchmark functions (Sect. \ref{sec3.5}) and on a highly non-linear FEA problem for verification (Sect. \ref{sec3.6}).
}

\subsection{Analytical test functions for optimization}\label{sec3.2}
{To test the performances of the investigated sampling methods, six two-dimensional benchmark tests are investigated. The formulae of the functions are listed below.
\\}

{Shubert 2D}
\begin{flalign}
    f(\bm{x})=& \prod_{i=1}^{2}{\left(\sum_{j=1}^5{ cos((j+1)x_i+j)}\right)} \:\:\:\:\:\:\:\:\:\ x_{1,2} \in [-2,2]&&
\label{eq12}
\end{flalign}

{Ackley 2D}
\begin{flalign}
    f(\bm{x}) = -20 \cdot e^{-0.2\sqrt{\frac{1}{2}\sum_{i=1}^{2}x_i^2}}-e^{\frac{1}{2}\sum_{i=1}^{2}cos(2 \pi \cdot x_i)}+ 20 + e \:\:\:\:\:\:\:\:\:\ x_{1,2} \in [-5,5]&&
\label{eq13}
\end{flalign}

{Rosenbrock 2D}
\begin{flalign}
    f(\bm{x})=\sum_{i=1}^{1}[100 (x_{i+1} - x_i^2)^ 2 + (x_i-1)^2] \:\:\:\:\:\:\:\:\:\ x_{1,2} \in [-2,2]&&
\label{eq14}
\end{flalign}

{Michalewicz 2D}
\begin{flalign}
    f(\bm x) = - \sum_{i=1}^{2}sin(x_i)sin^{20}\left(\frac{ix_i^2}{\pi}\right) \:\:\:\:\:\:\:\:\:\ x_{1,2} \in [0,4]&&
\label{eq15}
\end{flalign}

{Sphere 2D}
\begin{flalign}
    f(\bm{x})=\sum_{i=1}^{2} x_i^{2} \:\:\:\:\:\:\:\:\:\ x_{1,2} \in [-5,5]&&
\label{eq16}
\end{flalign}

{Zakharov 2D}
\begin{flalign}
    f(\bm{x})=\sum_{i=1}^2 x_i^{2}+(\sum_{i=1}^2 0.5ix_i)^2 + (\sum_{i=1}^2 0.5ix_i)^4 \:\:\:\:\:\:\:\:\:\ x_{1,2} \in [-10,10]&&
\label{eq17}
\end{flalign}

{The analytical functions listed above are often used as benchmark functions for optimization problems. According to the classification proposed by \citet{Winston1992}, these functions can be grouped in terms of features, such as modality, valleys, separability, and dimensionality. These features are shortly summarized below.}

\begin{itemize}
	\item \textbf{Modality}: A function with more than one local optimum is called multimodal, otherwise it is called unimodal.
	\item \textbf{Valleys}: Functions where a narrow area of little change is surrounded by regions of steep descent.
	\item \textbf{Separability}: This is a measure of the difficulty in optimizing a given benchmark function. In the literature different definitions of separability are given \citep{SALOMON1996263}, \citep{Boyer2011}. Generally speaking, functions that present inter-relations between variables are non-separable.
    \item \textbf{Dimensionality}: Number of dimensions (or variables) of a given function.
\end{itemize}

{
The choice of the presented six benchmark functions (Group 1) was made based on the features outlined in the work of \citet{Jamil2013ALS}, \citet{Adorio2005MVFM}, and \citet{Molga2005}.
To begin with, the Shubert function is a multimodal uniformly waving function and has several local minima. Similarly, the Ackley function is also a multimodal function with many local minima but is not separable. Moreover, the Ackley function is perfectly symmetric to each of its variables. This function is highly nonlinear and is characterized by a nearly flat outer region and a large hole at the center leading to its global minimum. Rosenbrock function is a unimodal, non-separable function, which has a single, nearly flat valley at the center. It is often used to test gradient-based optimization algorithms because of the flatness of the valley. The Michalewicz function has a number of local minima equal to the dimensionality of the function itself. The area containing the minima is very small compared to the entire search space. This means that this function has very steep ridges comparable to mathematical discontinuities. On the contrary, the Sphere function is separable, convex, and unimodal. Finally, the Zakharov function is a plate-shaped function without local minima except the global one. This function is multimodal and non-separable.
}

\begin{figure}[!htb]
\minipage{1.0\textwidth}
\centering
  \includegraphics[width=0.95\columnwidth]{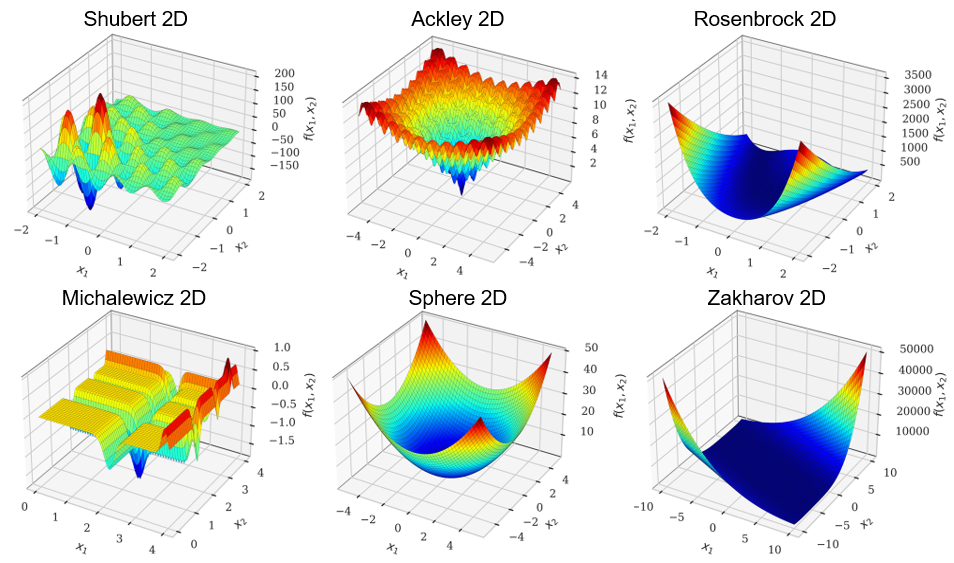}
  \caption{Two-dimensional Shubert, Ackley, Rosenbrock, Michalewicz, Sphere, and Zakharov function.}
  \label{fig_08}
\endminipage\hfill
\end{figure}

{
Because of the “curse of dimensionality”, most of the metamodeling techniques and optimization algorithms are generally affected by the increasing dimensionality of a problem. According to \citet{Winston1992} and \citet{Yao2003}, as the number of dimensions increases, the search space increases exponentially. This could be a relevant barrier for some sampling techniques as well.
\\
For the above reasons, it is worth investigating the performance of the sampling algorithms in problems with different dimensionality. The Ackley, Rosenbrock, and Sphere functions are investigated also in a 5-, 10- and 30-dimensional design space (Group 2). The generic d-dimensional formulae are given below.
}

{Ackley d-D}
\begin{flalign}
    f(\bm{x}) = -20 \cdot e^{-0.2\sqrt{\frac{1}{d}\sum_{i=1}^{d}x_i^2}}-e^{\frac{1}{d}\sum_{i=1}^{d}cos(2 \pi \cdot x_i)}+ 20 + e \:\:\:\:\:\:\:\:\:\ x_{1,...,d} \in [-5,5]&&
\label{eq18}
\end{flalign}

{Rosenbrock d-D}
\begin{flalign}
    f(\bm{x})=\sum_{i=1}^{d-1}[100 (x_{i+1} - x_i^2)^ 2 + (x_i-1)^2] \:\:\:\:\:\:\:\:\:\ x_{1,...,d} \in [-2,2]&&
\label{eq19}
\end{flalign}

{Sphere d-D}
\begin{flalign}
    f(\bm{x})=\sum_{i=1}^{d} x_i^{2} \:\:\:\:\:\:\:\:\:\ x_{1,...,d} \in [-5,5]&&
\label{eq20}
\end{flalign}

{A comprehensive listing of the benchmark functions considered is summarized in the table below.}

\begin{table}
	\caption{Summary of benchmark functions}
	\centering
	\begin{tabular}{lllll}
		\toprule
		% \multicolumn{2}{c}{Part}                   \\
		% \cmidrule(r){1-2}
		Group  & Benchmark function  & Domain & Number of variables & Description \\
		\midrule
		Group 1 & Shubert 2D  & $[-2,2]^2$ & 2 & Multimodal, high nonlinearity \\
		          & Ackley 2D   & $[-5,5]^2$ & 2 & Multimodal, high nonlinearity, symmetric  \\
		          & Rosenbrock 2D   & $[-2,2]^2$ & 2 & Unimodal, valley-shaped  \\
		          & Sphere 2D   & $[-5,5]^2$ & 2 & Unimodal, bowl-Shaped, symmetric  \\
		          & Zakharov 2D   & $[-10,10]^2$ & 2 & Unimodal, plate-Shaped  \\
		  Group 2 & Ackley 5D   & $[-5,5]^5$ & 5 & Medium dimensionality  \\
		          & Rosenbrock 5D   & $[-2,2]^5$ & 5 & Medium dimensionality  \\
		          & Sphere 5D   & $[-5,5]^5$ & 5 & Medium dimensionality  \\
		          & Ackley 10D   & $[-5,5]^{10}$ & 10 & High dimensionality  \\
		          & Rosenbrock 10D   & $[-2,2]^{10}$ & 10 & High dimensionality  \\
		          & Sphere 10D   & $[-5,5]^{10}$ & 10 & High dimensionality  \\
		        & Ackley 30D   & $[-5,5]^{30}$ & 30 & High dimensionality  \\
		          & Rosenbrock 30D   & $[-2,2]^{30}$ & 30 & High dimensionality  \\
		          & Sphere 30D   & $[-5,5]^{30}$ & 30 & High dimensionality  \\
		\bottomrule
	\end{tabular}
	\label{tab01}
\end{table}

\subsection{Test scheme}\label{sec3.3}
\begin{figure}[!htb]
\minipage{1.0\textwidth}
\centering
  \includegraphics[]{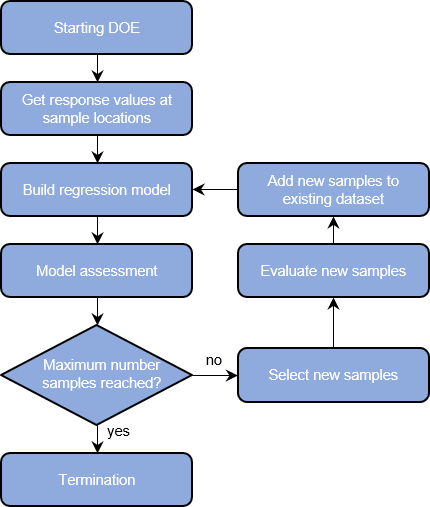}
  \caption{Flow chart for the numerical evaluation of sequential sampling methods.}
  \label{fig_09}
\endminipage\hfill
\end{figure}

{
As a pre-processing step, the initial domain is linearly scaled to the hypercube $[0,1]^d$. As for adaptive methods, a space-filling LHS is used as the starting dataset, which is subsequently refined until the maximum number of allowed samples has been reached. This choice is aimed at having a fair comparison with one-stage methods. Nevertheless, for optimal performance, it would be ideal to start with a pre-optimized LHD whenever applicable. At each iteration, a metamodel is trained on the available samples. To prove that the performances of the sampling strategies are independent of metamodeling techniques, two regression models are investigated: Kriging or Gaussian Process (GP) and Support Vector Regression (SVR). Kriging is employed with a Rational Quadratic kernel while the SVR method with a Radial Basis Function kernel. A detailed description of GP and SVR can be found in \citep{duvenaud-thesis-2014}, \citep{Rasmussen2006} and \citep{Brereton2010SupportVM}, \citep{Gunn1998SupportVM} respectively. A schematic overview of the procedure used for numerical evaluation of sequential methods is represented by the flow chart in Figure \ref{fig_09}.
}

{
Concerning the benchmark functions, in order to evaluate the accuracy of the regression model, the Root Mean Square Error (RMSE), defined in Eq. \ref{eq21}, is employed.
}

\begin{equation}
    \label{eq21}
    RMSE = \sqrt{\frac{1}{t}\sum_{i=1}^{t}(f(\bm{x}_i)-\hat{f}(\bm{x}_i))^2}
\end{equation}

{
$f$ is the true function, $\hat{f}$ is a metamodel built on the current samples and $t$ represents the number of randomly spread test points (here $5000 \cdot d$) in the domain. Although this measure is quite accurate, it is not applicable for expensive black box functions such as FEM simulations.\\
The RMSE results are compared with the RMSE of a one-stage LHS, which is one of the most widely used sampling strategies for dealing with expensive black box optimizations. Since pre-optimized LHDs are not available for every combination of dimensions and sampling points, the best space-filling design from $1000 \cdot d$ Monte Carlo-based LHS design is considered. Non-collapsing properties are not expressly optimized, but keep in mind that the LHDs already benefit from good projective properties by definition. The results are marked as sf-LHS. This sampling strategy is also the same one that is used to generate the initial dataset for adaptive methods.\\
Python 3.8 is the programming language used to get the results.
}

\subsection{Results}\label{sec3.4}
{The resulting mathematical properties are summarized in Figure \ref{fig_10}. Regarding MIPT, the parameter $\alpha$ is automatically tuned using the approach proposed in Sect. \ref{sec2.5.1}. The following design spaces were investigated:}

\begin{itemize}
	\item 2-dimensional design space with 20 samples 
	\item 5-dimensional design space with 50 samples 
	\item 10-dimensional design space with 100 samples
\end{itemize}

\begin{figure}[!htb]
\minipage{1.0\textwidth}
\centering
  \includegraphics[width=0.7\columnwidth]{pics/fig_10.png}
  \caption{Intersite and projected distance in 2-, 5-, and 10-dimensional design spaces with 20,50, and 100 samples respectively.}
  \label{fig_10}
\endminipage\hfill
\end{figure}

{It is worth mentioning that the bar plots in Figure \ref{fig_10} are solely dependent on the dimensionality of the domain and the size of the dataset. They are therefore unaffected by the problem investigated.\\ 
The pre-optimized LHD outperforms any other algorithms in intersite- and projected distance, which is expected, since these designs have been intensively optimized, and the number of samples is known a priori. However, this relevant information is often not available up-front. On the other hand, adaptive methods are specifically designed to be unrelated to the size of the DoE. If this information is available, pre-optimized LHS should be the first choice for sampling.\\
Mc-intersite-proj-th confirms its potential as a powerful sequential algorithm. The investigation shows that it outperforms mc-intersite-proj with respect to the selected mathematical criteria. This also suggests that the automatic adjustment of the $\alpha$ parameter appears to work appropriately as the dimensionality of the problem changes.\\
As for the new methods proposed in Sect. \ref{sec2.5}, they show on average great space-filling properties. The projective properties are, however, relatively poor. This makes sense since the projective properties have not been explicitly optimized. Note that the proposed designs are very close to LHDs and have very similar desirable mathematical properties as a result. Furthermore, these methods are on average comparable to, if not slightly preferable to, one-stage LHDs.\\
As far as low-discrepancy sequences are concerned, the Sobol method shows very good projective properties. Among the samplings investigated here, it has the second-best non-collapsing properties by a clear margin, regardless of the dimensionality of the problem and the number of samples. However, especially for a problem with a small number of samples, its space-filling properties are relatively poor, even when compared to Halton sampling.\\
Randomly generated LHSs showed the poorest performance. Apart from the slight improvement in terms of intersite distance offered by the improved Monte Carlo-based version (sf-LHS), the random distribution of samples in the grid intervals makes the non-collapsing properties mediocre, which is quite unusual for Latin Hypercube Designs. Moreover, this can be seen as unfavorable in practical applications such as crashworthiness and optimizations where certain input variables have varying effects on response results. 
}

\subsection{Application on the benchmark functions}\label{sec3.5}
{Due to their simple and efficient generation, LHSs are (and will likely remain) the most widely used one-shot sampling approaches. Therefore, a deeper study is conducted on benchmark functions (presented in Sect. \ref{sec3.2}) to compare the sf-LHS with the most promising adaptive sampling method, i.e. MIPT. To provide a basis for comparison, the new FpPLHS and MqPLHS methods are also included. They could offer a viable alternative when it is expressly desired to have a quasi-LHD. Also, since they are independent of any parameter tuning, there is interest in seeing their performance for medium to high dimensionality.
}
\begin{figure}[!htb]
\minipage{1.0\textwidth}
\centering
  \includegraphics[width=1.0\columnwidth]{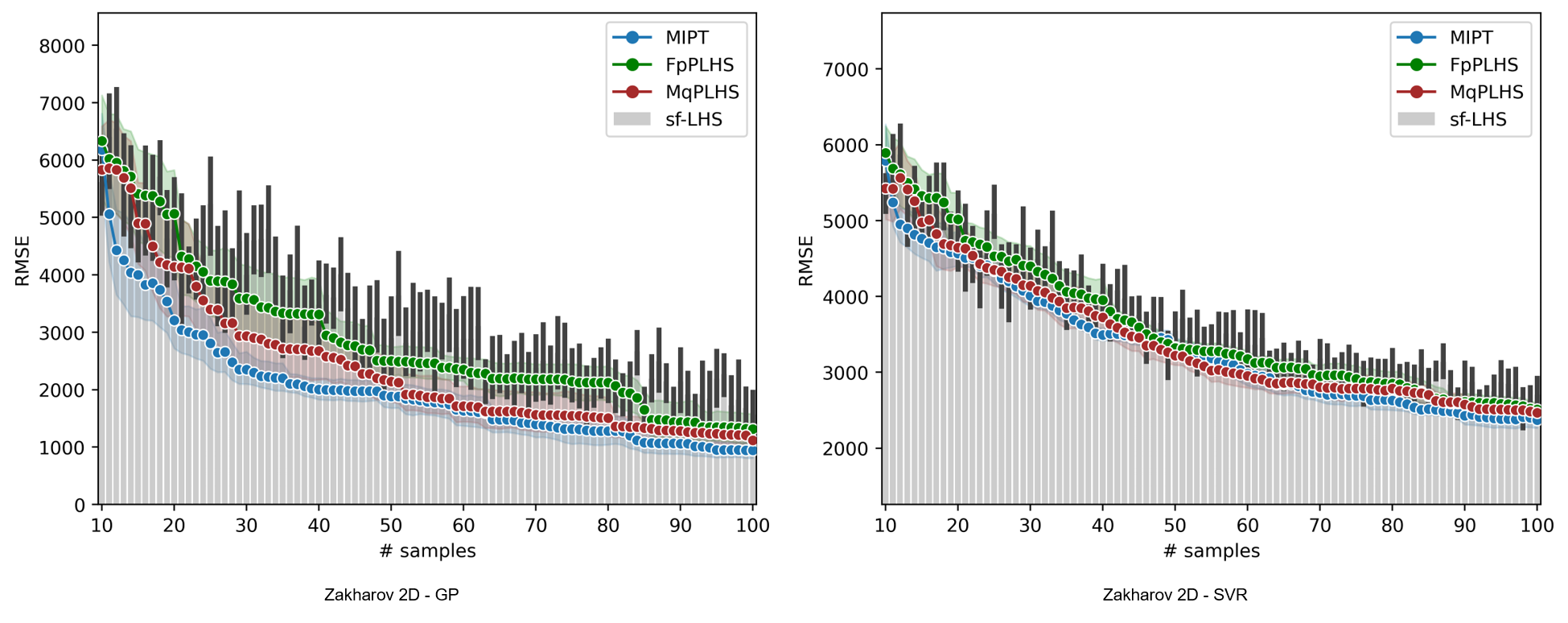}
  \caption{RMSE error of adaptive methods MIPT, FpPLHS, MqPLHS (solid lines) and sf-LHS (bar plot) over sample size for the 2D Zakharov function with Gaussian Process (left) and SVR (right) metamodels.}
  \label{fig_11}
\endminipage\hfill
\end{figure}

{
By way of example, in Figure 11, the RMSE curves obtained for the 2D Zakharov function are shown. The two metamodeling techniques, GP and SVR are depicted on the left and the right respectively. After some initial adjustment iterations, the RMSE error of the sequential methods decreases gradually, regardless of the adaptive method and metamodeling technique chosen. All three adaptive methods outperform the sf-LHS method in terms of both mean error (smaller) and interval bands (much narrower). The MIPT algorithm delivers the best performance by far in this benchmark function, further corroborating the tests observed in Sect. \ref{sec3.4}. Even FpPLHS and MqPLHS yield encouraging results in terms of mean error, but with slightly wider confidence intervals than MIPT. The average RMSE value of the FpPLHS method also shows a sporadic trend that seems to resemble the frequency of the objective function observed in Figure \ref{fig_07}. By taking the worst result among the adaptive methods, an improvement of about 19\% and 9\% at the final stage is still guaranteed for GP and SVR respectively.
\\\\
This result underlines the potential of sequential methods, especially since the final number of samples does not have to be known in advance. Additionally, the shape of the convergence curve shows a flattening characteristic on average after 65-70 iterations. Also, the standard deviation value remains relatively small. These are ideal properties for the application of a cut-off convergence criterion.
\\\\
A further example of a challenging benchmark function where the benefits of adaptive methods are less obvious is given by the 2D Ackley function.
}

\begin{figure}[!htb]
\minipage{1.0\textwidth}
\centering
  \includegraphics[width=1.0\columnwidth]{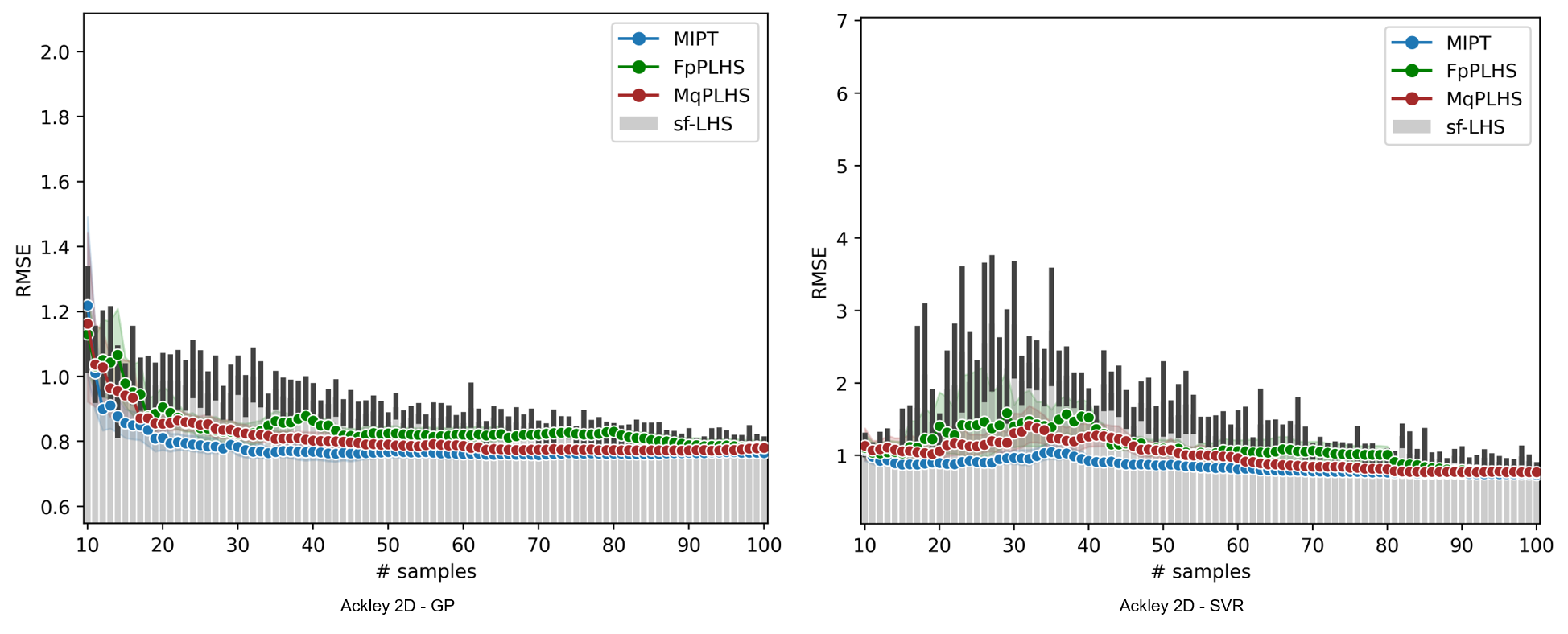}
  \caption{RMSE error of adaptive methods MIPT, FpPLHS, MqPLHS (solid lines), and sf-LHS (bar plot) over sample size for the 2D Ackley function with Gaussian Process (left) and SVR (right) metamodels.}
  \label{fig_12}
\endminipage\hfill
\end{figure}

\begin{figure}[!htb]
\minipage{1.0\textwidth}
\centering
  \includegraphics[width=1.0\columnwidth]{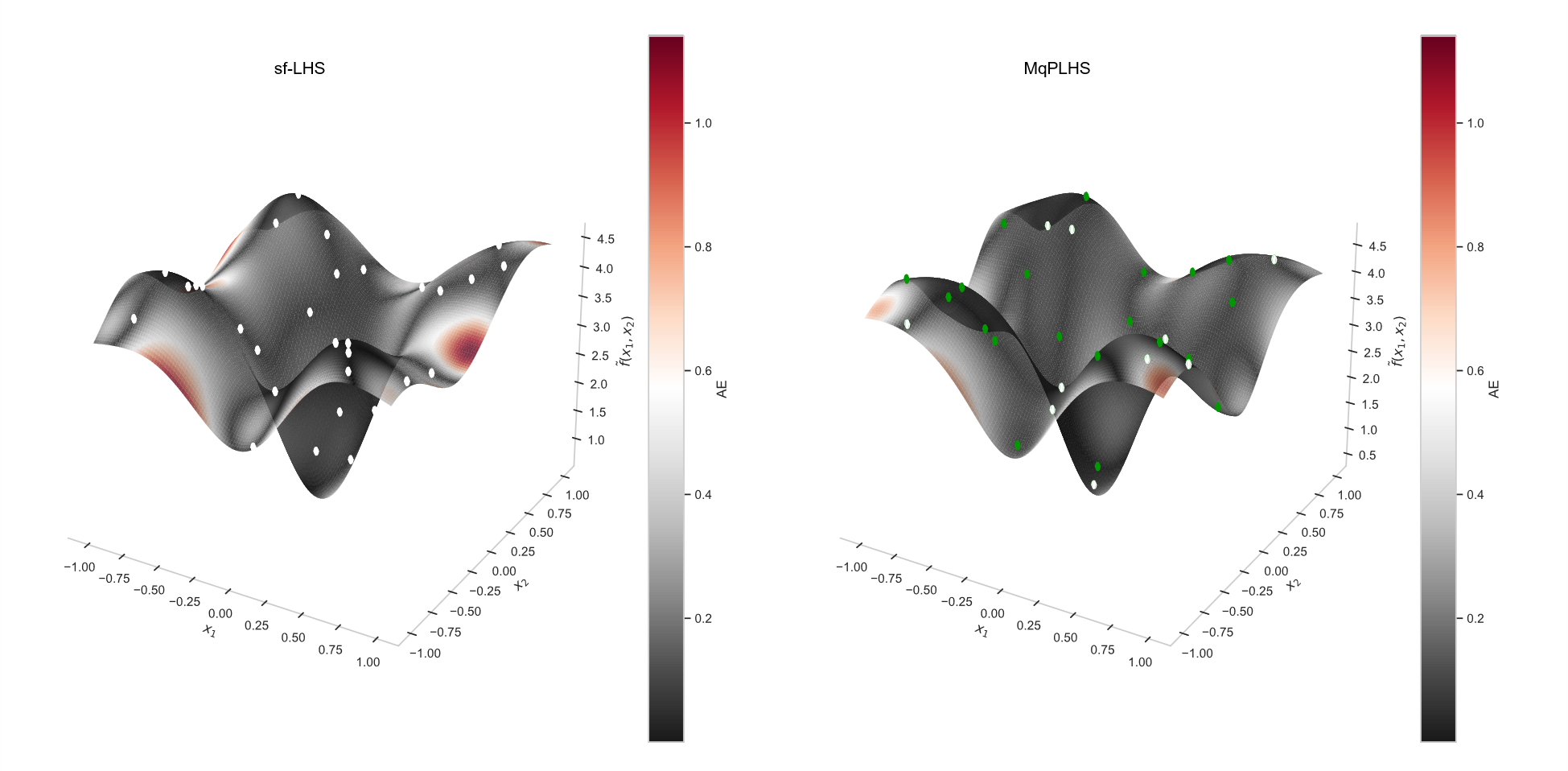}
  \caption{Comparison of the absolute error (AE) of the Ackley function with 30 samples using sf-LHS (on the left) and MqPLHS (on the right) sampling methods.}
  \label{fig_13}
\endminipage\hfill
\end{figure}

{
As shown in Figure \ref{fig_12}, albeit to a smaller extent, adaptive methods outperform the sf-LHS method in each way, particularly in terms of confidence intervals. The convergence rate is not as clear as in the previous function. This is most likely due to the difficulty of the metamodels in matching the numerous local minima of this function with so few samples available. In addition, the notable difference between the two metamodeling techniques is probably due to the greater ease of Kriging in approximating the high non-linearities. Because this function is perfectly symmetric and uniformly wavy, the differences between the sampling strategies are less evident. Theoretically, space-filling properties are preferred over projected properties to achieve adequate accuracy results within this benchmark function.
\\\\
To convince the reader of the better accuracy achieved by the adaptive methods in this function, a further visualization of the metamodels is provided in Figure \ref{fig_13}. The GP metamodels for an sf-LHS and MqPLHS design are shown in a $[-1,1]^2$  with 30 samples. White dots represent initial samples that were used with a one-shot method (here sf-LHS). The green points instead are the result of a refinement of successive iterations of the MqPLHS method. The colormap used shows the trend of the absolute error compared to the real function: the more a region is colored red, the greater the error in that part of the function.
}

\begin{figure}[!htb]
\minipage{1.0\textwidth}
\centering
  \includegraphics[width=1.0\columnwidth]{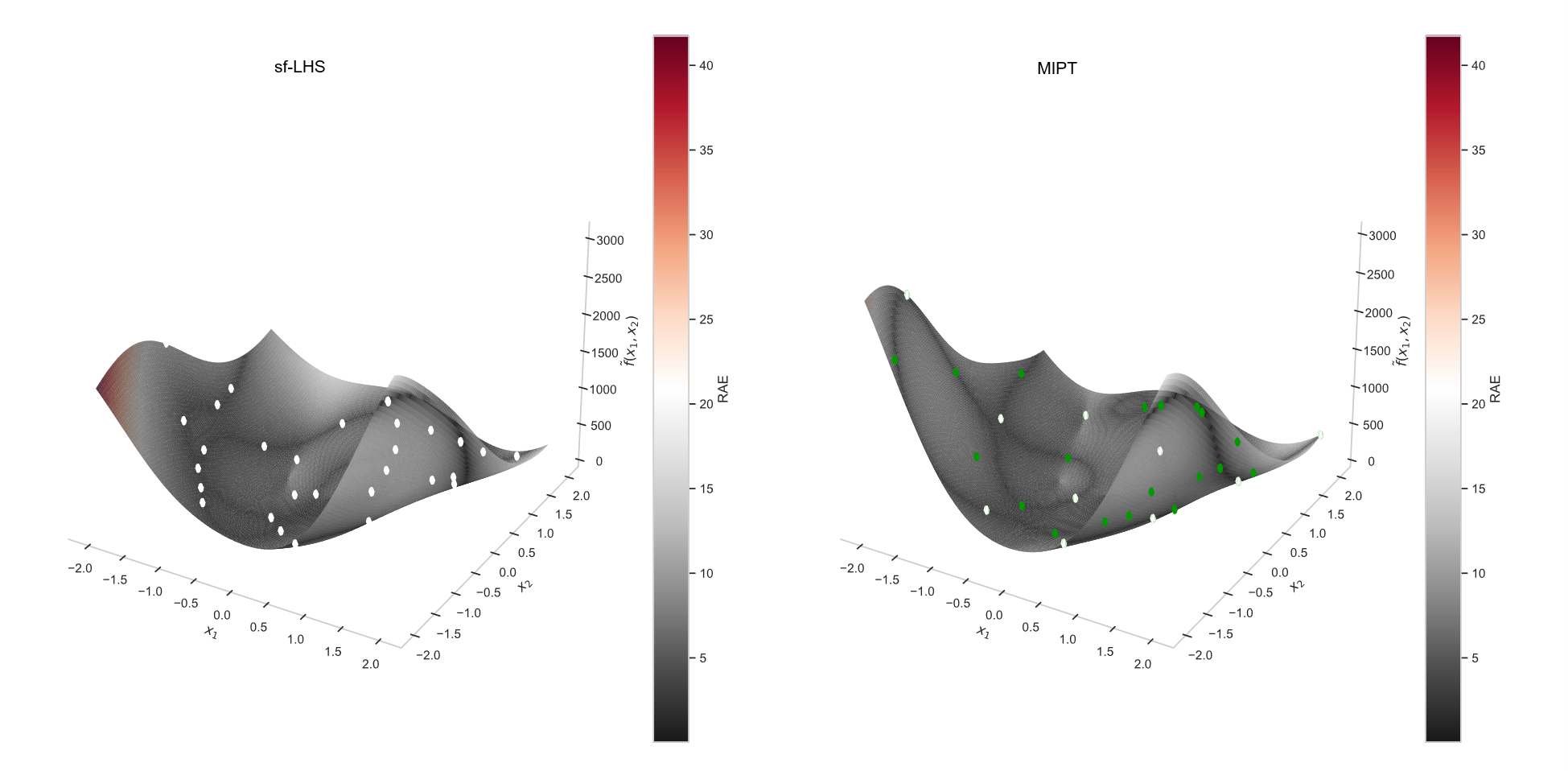}
  \caption{Comparison of the absolute error (AE) of the Rosenbrock function with 30 samples using sf-LHS (on the left) and MIPT (on the right) sampling methods.}
  \label{fig_14}
\endminipage\hfill
\end{figure}

{
The same visualization can be equally appreciated for other benchmark functions, other sampling methods, error metrics, and metamodeling techniques. For example, Figure \ref{fig_14} shows the SVR-approximation of the Rosenbrock function with 30 samples in the $[-2,2]^2$ domain. Sf-LHS and MIPT are compared. To ideally scale the values, here the colormap accounts for the square root of the absolute error (RAE).
}

\begin{figure}[!htb]
\minipage{1.0\textwidth}
\centering
  \includegraphics[width=1.0\columnwidth]{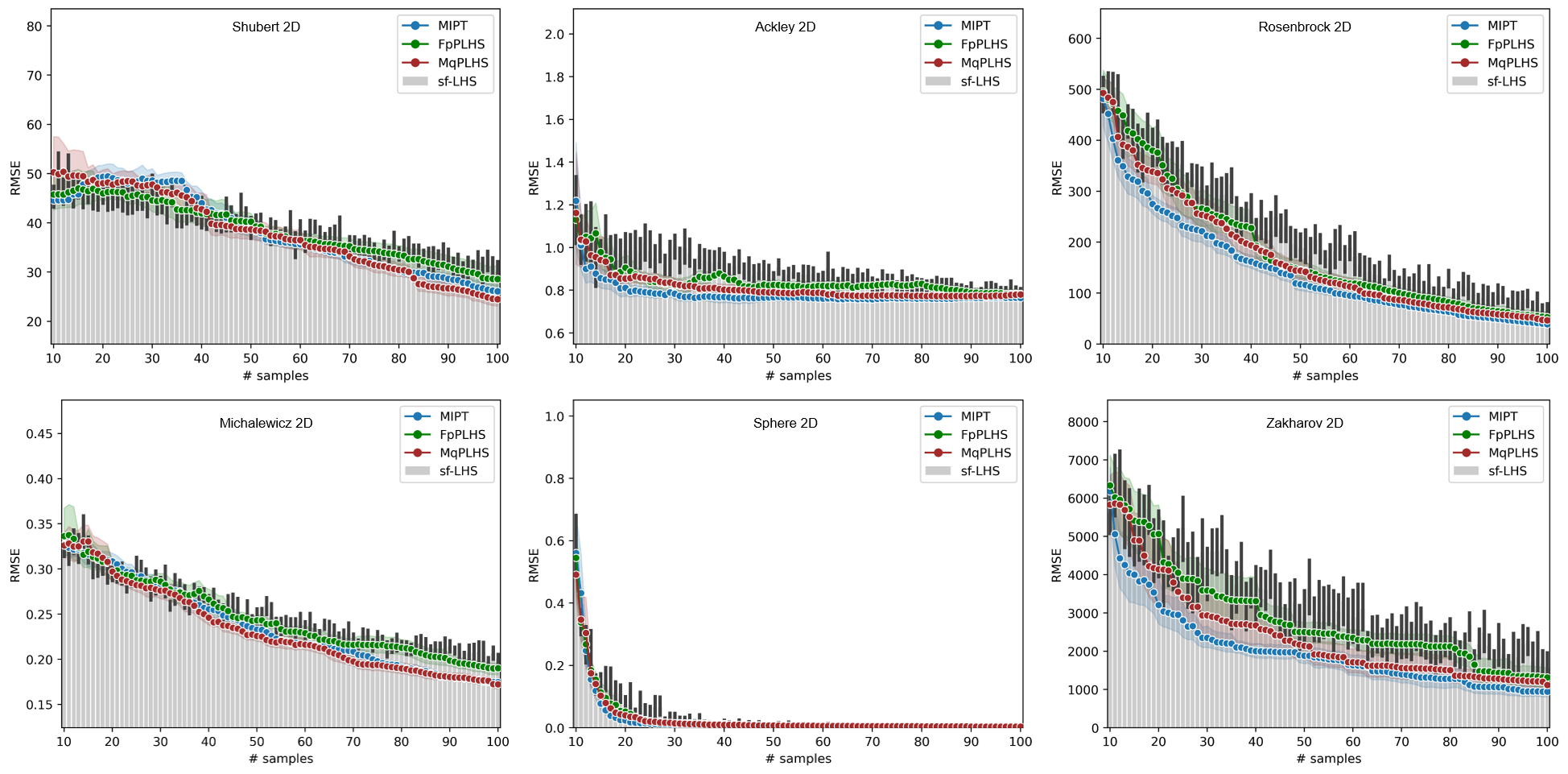}
  \caption{2-dimensional benchmark functions modeled with GP}
  \label{fig_15}
\endminipage\hfill
\end{figure}

\begin{figure}[!htb]
\minipage{1.0\textwidth}
\centering
  \includegraphics[width=1.0\columnwidth]{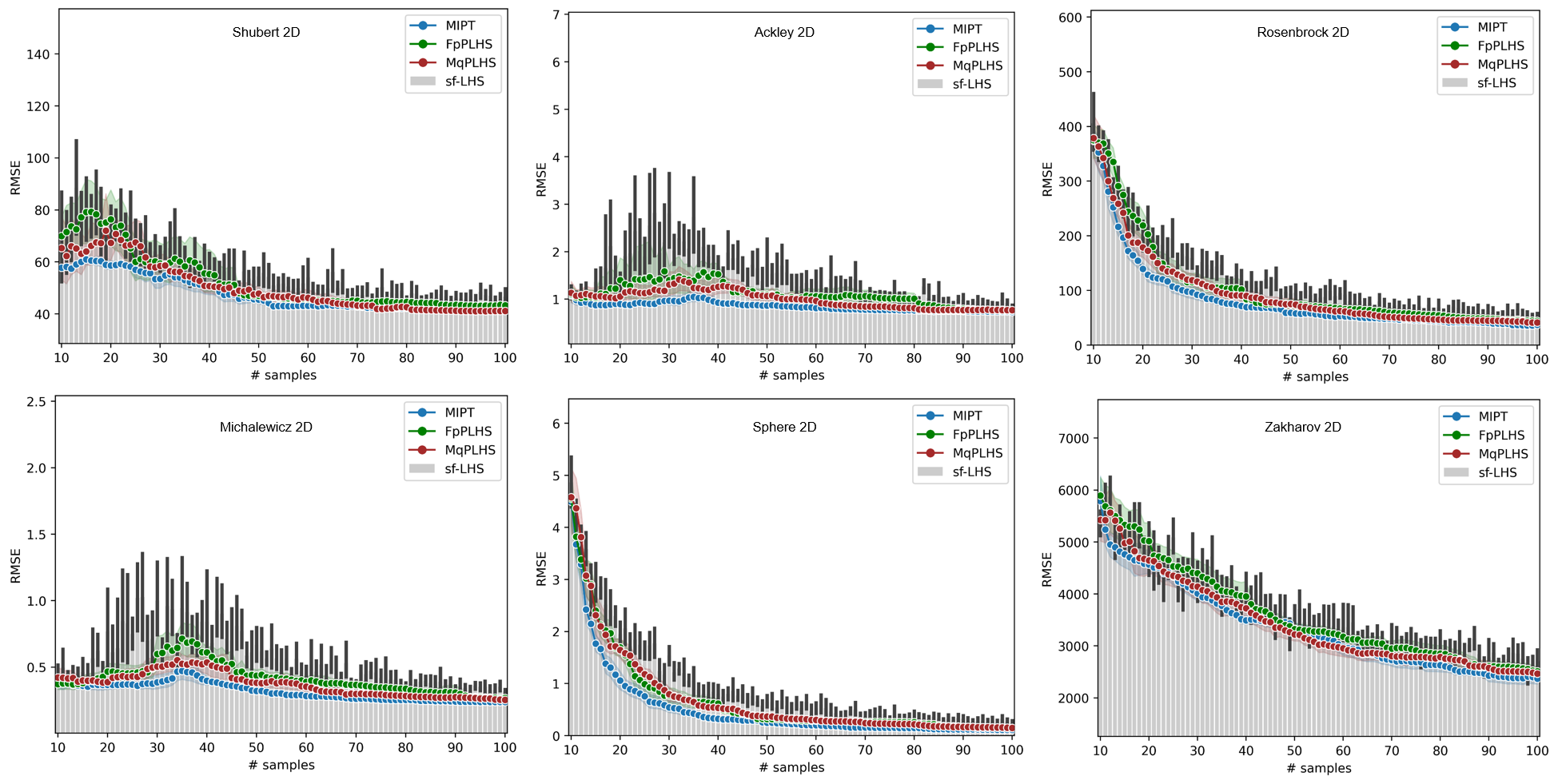}
  \caption{2-dimensional benchmark functions modeled with SVR}
  \label{fig_16}
\endminipage\hfill
\end{figure}

\begin{figure}[!htb]
\minipage{1.0\textwidth}
\centering
  \includegraphics[width=1.0\columnwidth]{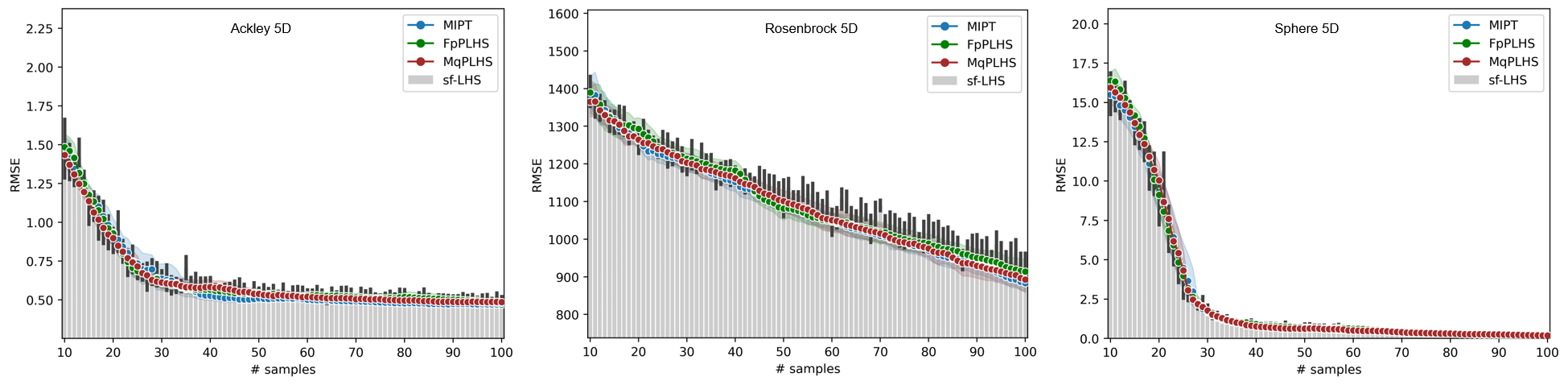}
  \caption{5-dimensional benchmark functions modeled with GP}
  \label{fig_17}
\endminipage\hfill
\end{figure}

\begin{figure}[!htb]
\minipage{1.0\textwidth}
\centering
  \includegraphics[width=1.0\columnwidth]{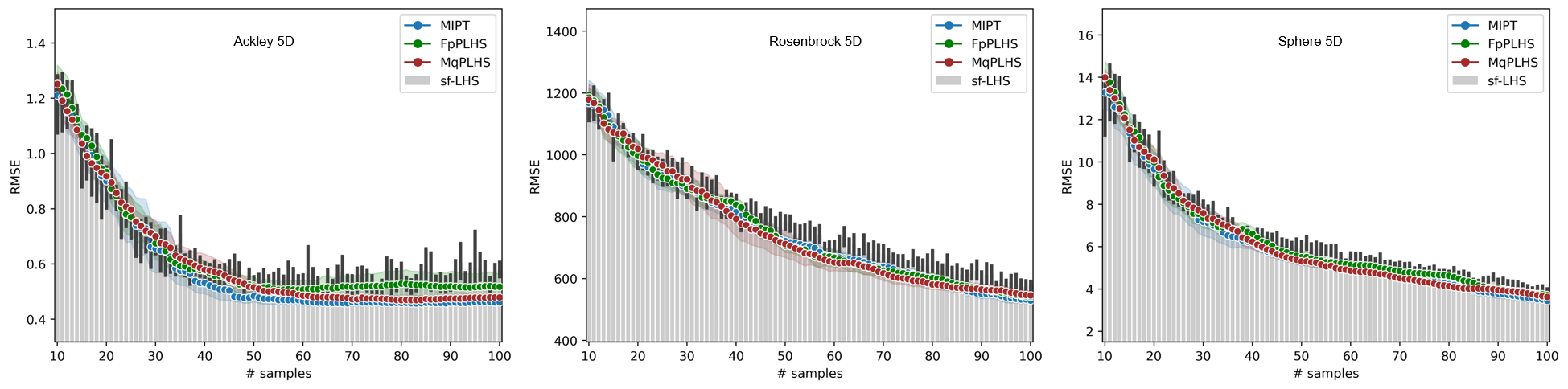}
  \caption{5-dimensional benchmark functions modeled with SVR}
  \label{fig_18}
\endminipage\hfill
\end{figure}

\begin{figure}[!htb]
\minipage{1.0\textwidth}
\centering
  \includegraphics[width=1.0\columnwidth]{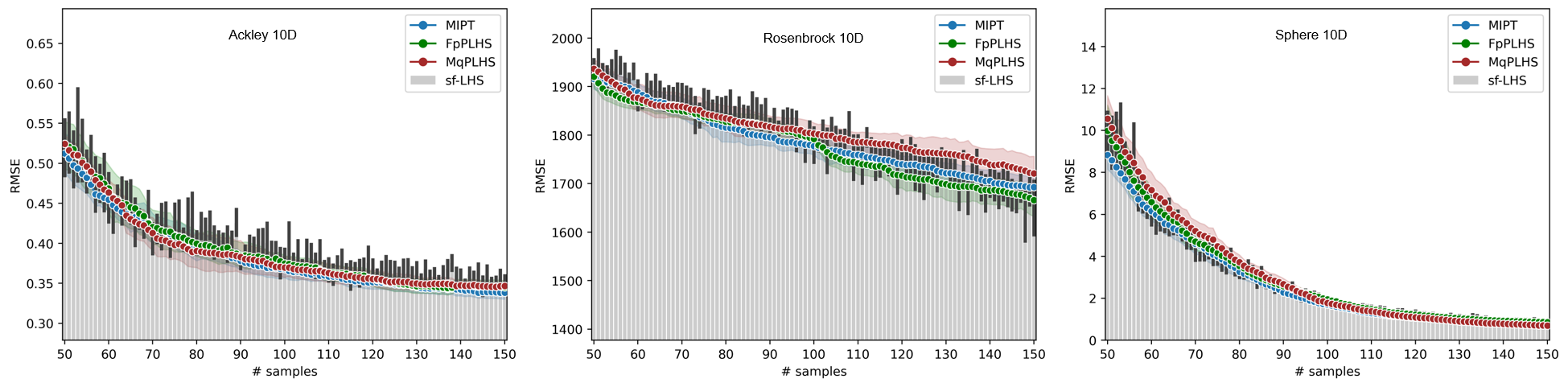}
  \caption{10-dimensional benchmark functions modeled with GP}
  \label{fig_19}
\endminipage\hfill
\end{figure}

\begin{figure}[!htb]
\minipage{1.0\textwidth}
\centering
  \includegraphics[width=1.0\columnwidth]{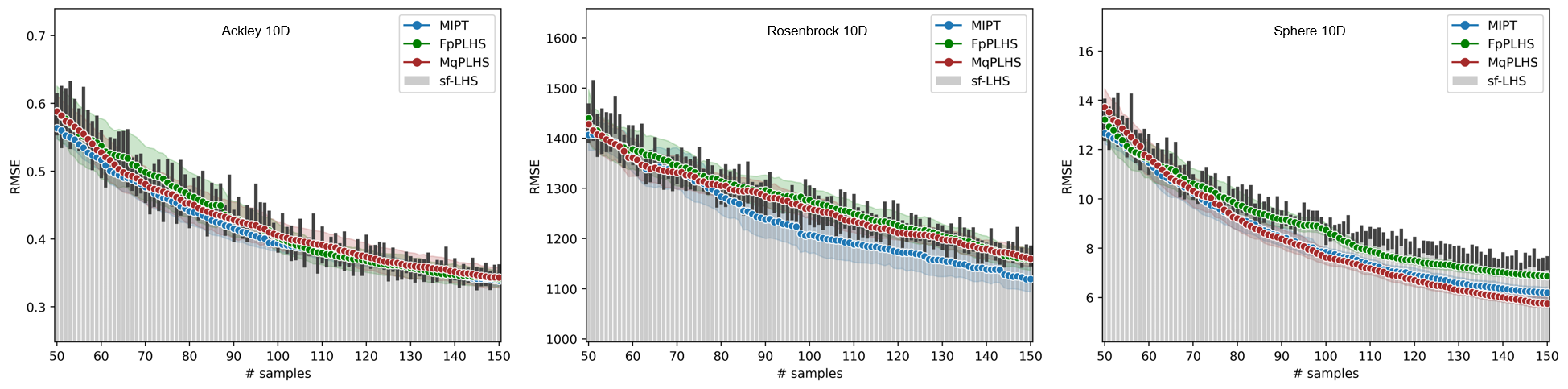}
  \caption{10-dimensional benchmark functions modeled with SVR}
  \label{fig_20}
\endminipage\hfill
\end{figure}

\begin{figure}[!htb]
\minipage{1.0\textwidth}
\centering
  \includegraphics[width=1.0\columnwidth]{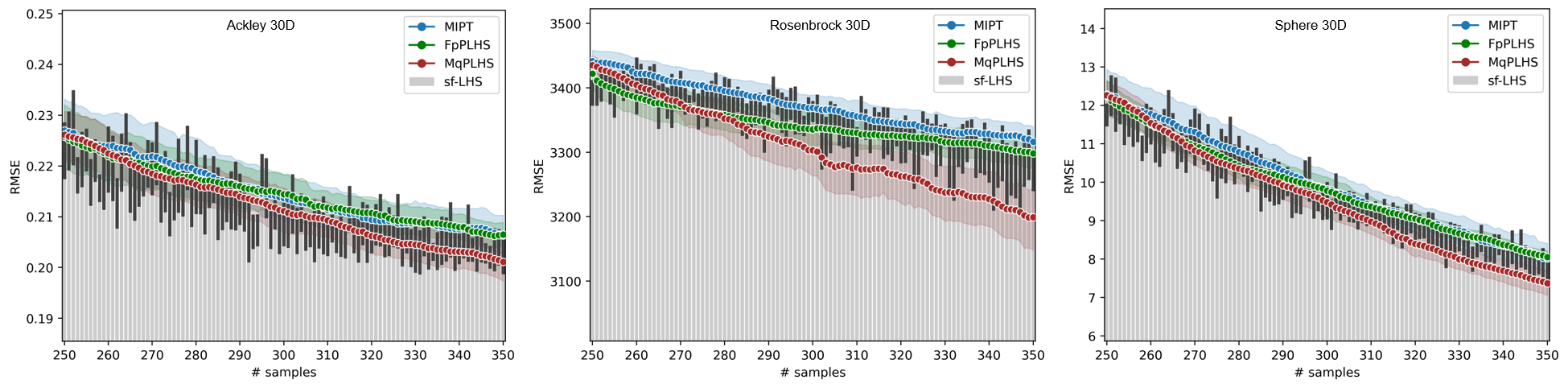}
  \caption{30-dimensional benchmark functions modeled with GP}
  \label{fig_21}
\endminipage\hfill
\end{figure}

\begin{figure}[!htb]
\minipage{1.0\textwidth}
\centering
  \includegraphics[width=1.0\columnwidth]{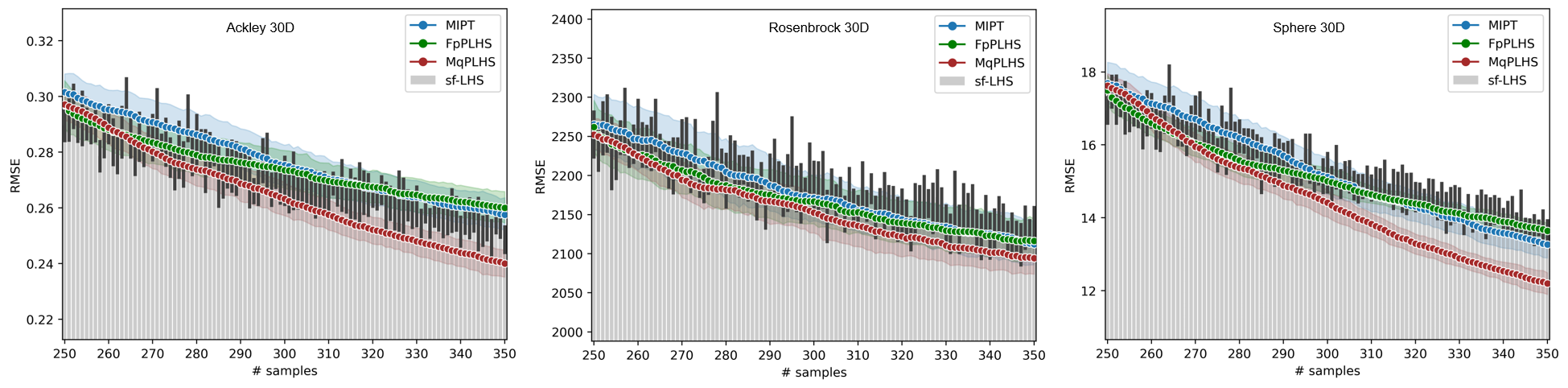}
  \caption{30-dimensional benchmark functions modeled with SVR}
  \label{fig_22}
\endminipage\hfill
\end{figure}

{Since it is not feasible to draw general conclusions from only a few functions, all benchmark functions presented in Sect. \ref{sec3.2} are tested with both GP and SVR regression methods. The results are shown in Figures \ref{fig_15}-\ref{fig_22}.
A first remark concerns the two metamodeling techniques used. Regardless of the function used and the dimensionality of the problem, the results seem consistent between the two regression methods. This suggests that adaptive sampling methods have some benefits over a one-stage approach such as sf-LHS independently of the metamodeling technique. Therefore, the comments that follow will apply to both GP and SVR. \\
Although with different performances, all three proposed adaptive methods clearly outperform space-filling LHS in two dimensions. MIPT confirms its superiority in terms of convergence speed, RMSE values, and confidence intervals. This performance is very close to the one exhibited by MqPLHS. In particular, MqPLHS results are completely comparable and at times slightly superior in the Shubert and Michalewicz functions. As a bit of a surprise compared to the tests in Sect. \ref{sec3.4}, FpPLHS shows the poorest and most fluctuating performance of these three sampling methods. A final plateau is reached in every function, further confirming that a stopping algorithm would be ideal.\\
As the dimensionality of the problem grows, the difference between the adaptive methods and sf-LHS tapers off sharply, especially for the MIPT and FpPLHS methods. This is most likely due to the so-called "curse of dimensionality" effect that affects not only regression models \citep{Verleysen2005TheCO} but also sampling methods \citep{Forrester2008}. In this regard, the MqPLHS method seems to be the method that succeeds to cope with this issue the best. The difference in performance at the final stage (i.e., after 350 samples) in problems with 30 variables is remarkable (up to 14.3 \% improvement). Furthermore, with increasing dimensionality, the flatness of the observed plateau also decreases. This is indeed a factor to consider when implementing the stopping algorithm.
}

\subsection{Application on a crash-box optimization}\label{sec3.6}
{
A final verification test with respect to the sampling strategies under consideration is carried out on a generic structural component subject to highly non-linear deformations. A crash box, shaped as a square-based pyramidal frustum, is compressed by a rigid plane impacting in the vertical direction. As boundary conditions, the velocity of the rigid body is fixed to a constant value of 75 mm/s, and the main base of the pyramidal frustum is tied to a rigid plate. The impact between the crash box and the rigid plane is investigated for 0.07 seconds. The FEM model is composed of 2880 shell elements. The simulation is performed with the explicit solver LS-DYNA using 4 CPUs of an Intel Xeon W-2135 (8.25M Cache, 3.70 GHz) processor. The complete computational time for a single solver run is roughly 76 seconds.
}

\begin{figure}[!htb]
\minipage{1.0\textwidth}
\centering
  \includegraphics[width=0.8\columnwidth]{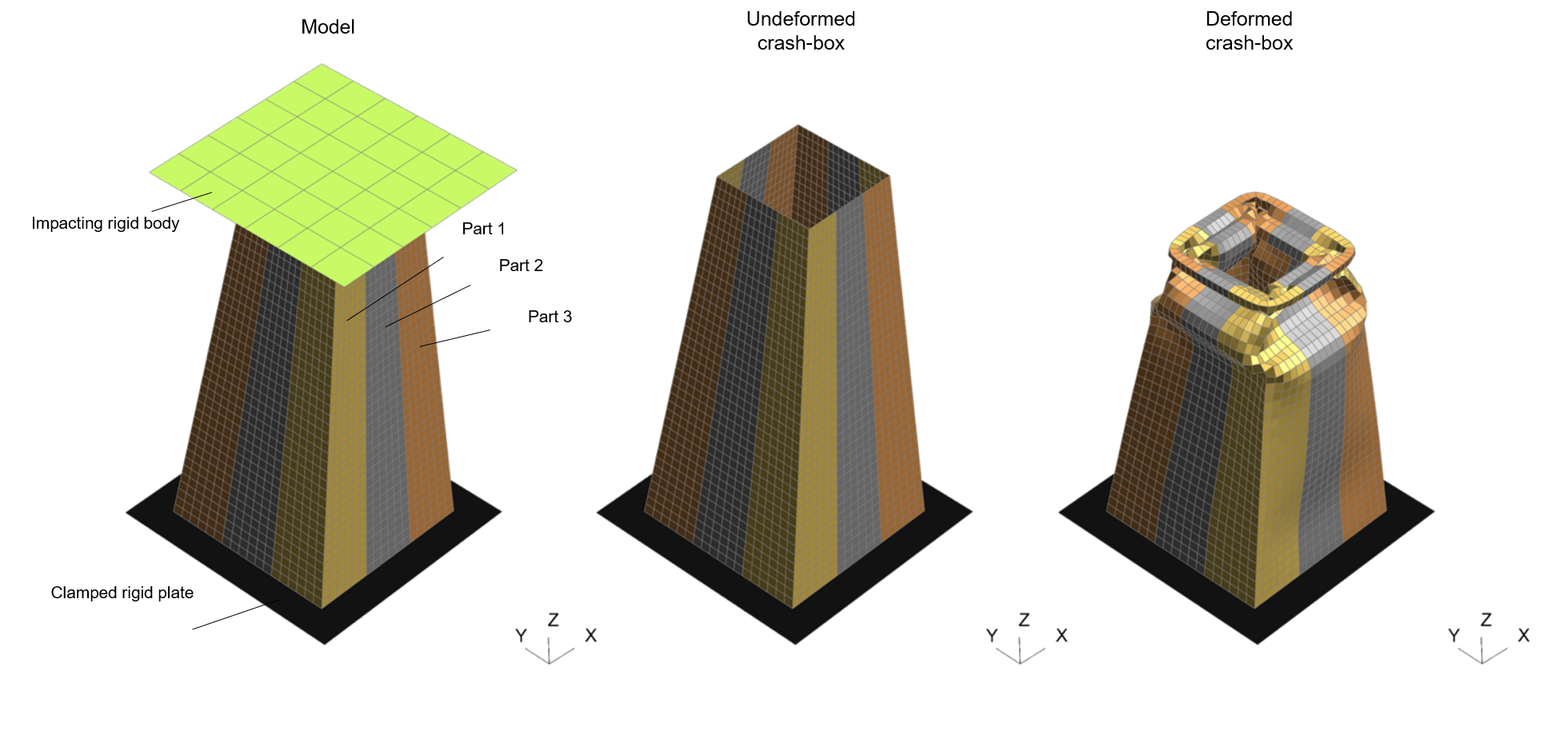}
  \caption{from left to right: undeformed crash box with impacting rigid body, undeformed crash box at $t=0.0$ $s$ and deformed crash box at $t=0.07$ $s$.}
  \label{fig_23}
\endminipage\hfill
\end{figure}

{As shown in Figure \ref{fig_23}, the crash box is vertically partitioned into three parts (colored in three different colors, respectively). The thickness of each part is used as the input variable. The mass of the structural component is the objective function, and the absorbed internal energy is set as a constraint function. Since optimization is beyond the scope of this paper, only the estimated model accuracy of these two functions is considered to evaluate the performances of the investigated sampling methods. Due to computational time restrictions, the application of each sampling strategy is repeated 5 times, and LHS is only computed with a step-size of 10 samples.
\\\\
Here, the error metric formula of Eq. \ref{eq21} is slightly modified by the employment of $k$-Fold Cross-Validation ($k=10$). In this way, there is no longer a need to evaluate $5000 \cdot d$ random points (barely feasible when it comes to FEM simulations). The samples of the dataset can be re-used through an appropriate train-test split strategy. The readjusted formula of this error metric, denoted by $RMSE_CV$, is shown in Eq. \ref{eq22}: 
}

\begin{equation}
    \label{eq22}
    RMSE_{CV} = \sqrt{\frac{1}{k}\sum_{i=1}^{k}(f(\bm{x}_i)-\hat{f}(\bm{x}_i))^2}
\end{equation}

\subsection{Discussion of the simulation results}\label{sec3.7}
{
For the crash box optimization, two response functions are required: the total mass of the component (objective function), and the internal energy at the final stage (constraint function) i.e., the absorbed energy during the deformation process. Since the mass of the component is linearly dependent on the thickness of the parts, the associated surrogate model can be accurately approximated. The total mass can be expressed by the following function:
}
\begin{equation}
    \label{eq23}
    m_{tot}=\rho_1 \cdot V_1+\rho_2\cdot V_2+\rho_3\cdot V_3 = \rho_1\cdot A_1\cdot t_1+\rho_2\cdot A_2\cdot t_2+\rho_3 \cdot A_3 \cdot t_3
\end{equation}

\begin{figure}[!htb]
\minipage{1.0\textwidth}
\centering
  \includegraphics[width=0.65\columnwidth]{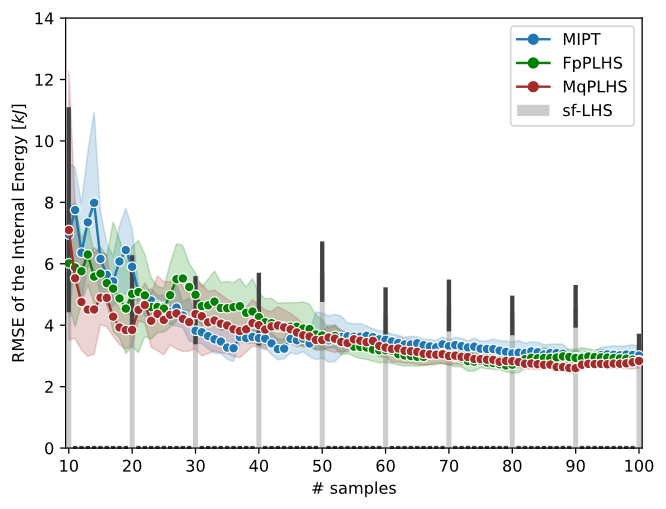}
  \caption{MRSE error of mc-intersite-proj-th (solid line) and sf-LHS (bar plot) for the crash box crushing test dependent on the sampling size.}
  \label{fig_24}
\endminipage\hfill
\end{figure}

{
where $\rho_i,V_i,A_i,t_i$ $for$ $i=1,2,3$ stands for density, volume, area, and thickness of the $i-th$ part, respectively. The areas of the parts and the densities of the material are constant. In this case, a linear surrogate model could deliver a more accurate approximation of the objective function than GP or SVR. Therefore, the focus is only on the internal energy response function. The comparison in Figure \ref{fig_24} shows that the adaptive sampling outperforms the one-shot approach, both in terms of standard deviation and average error. Even if, at the final stage, the improvement is limited, the adaptive sampling strategy reaches a plateau, which suggests that additional iterations will not lead to significant improvements. This obtained information can be fed into a convergence criterium to stop the Design of Experiments after a certain quality of the surrogate model is reached. This can avoid unnecessary solver runs, which reduces the computational costs, and can support the automated development process of numerical problems.  However, note that the cross-validation error curve is noisier in this case. This might make it challenging to apply an effective halting criterion for an early termination of the sampling process.
}

\section{Conclusions}\label{sec4}
{
Within the range of exploration-oriented sampling methods, sequential and one-stage approaches were investigated. The results confirm that sequential methods, under certain circumstances, provide an ideal sampling strategy, compared to the classic one-stage sampling methods, such as Latin Hypercubes. On the one hand, the results in terms of metamodel accuracy are generally improved. On the other hand, adaptive methods enable us to evaluate the quality of the surrogate model at each iteration step, which can be used for the application of a convergence criterion to stop the sampling process. By applying a convergence criterion, critical under- and oversampling phenomena can be avoided. 
\\\\
Regarding the three mathematical criteria selected for the evaluation of a sampling method (space-filling properties, non-collapsing properties, and granularity), the simplified version of the adaptive algorithm that aims to optimize both space-filling and non-collapsing properties (MIPT) provided optimal results overall. The empirical method proposed to automatically adjust the $\alpha$ parameter seems to cope well with dimensionality variations. This ensures that an optimal trade-off between space-filling - projective properties is achieved and that the algorithm always returns acceptable solutions.
\\\\
FpPLHS and MqPLHS were proposed with the aim of providing fine-grained approaches that could emulate an adaptive LHS method as closely as possible. 
\\\\
Whenever the number of samples should be known a priori, the best sampling strategy would be the pre-optimized LHS. It outperformed any other sampling method in terms of either projected- or intersite distance. However, they have been extensively optimized for several hours, and they are not available for certain design space dimensions and dataset sizes. This considerably limits their application but does not exclude them from being used as initial designs of sequential methods. On the contrary, classic LHS has shown low performance in the field of space-filling and only moderate non-collapsing properties. Unoptimized LHS should therefore be avoided.
\\\\
The verification tests on the benchmark functions partially confirm the expectations that arose from the criteria considered. More stable and more robust results are generally observed for sequential sampling strategies. Smooth curves featuring flattening plateaus allow the implementation of a convergence criterion to avoid expensive and unnecessary solver runs. Again here, especially for small and medium dimensional problems, MIPT showed great performance in terms of convergence rate, RMSE, and confidence interval.\\
MqPLHS has shown very promising results, especially in high-dimensionality problems. In functions with 30 variables, it seems to be the sampling method least affected by the "curse of dimensionality". An improvement of RMSE up to roughly 14\% was observed in the final stage compared to the other adaptive methods. The FpPLHS method, on the other hand, showed a relatively poor performance compared to the other two sequential methods, but still superior to the sf-LHS. Nonetheless, FpPLHS has a very simple implementation, is efficient in generating datasets for medium- and high-dimensionality problems ($d \geq 10$), and guarantees exact LHD at some stages, following the evolution of an exponential series.
\\\\
For verification, these methods are applied to the optimization of a crash box in the explicit simulation environment. The adaptive sampling methods show a good performance in terms of average error. The final plateau of the error curve obtained is the most interesting aspect of the adaptive sampling strategy. By applying a convergence criterion, a significant number of solver calls can be saved depending on the method considered. Especially in the field of crash simulations, where every iteration is usually computationally expensive, adaptive methods can significantly reduce the computational costs and avoid critical undersampling, which can make a subsequent optimization challenging. However, the cross-validation error curve is in this case noisier, which demands an appropriate stopping criterion that can deal with this potential issue.
\\\\
The following problems are the object of future research. First, a suitable and effective convergence criterion that can successfully deal with noisy error curves is also required. This halting algorithm will probably have to consider: the variation of the error throughout iterations (strongly influenced by the dimensionality of the problem), an upper limit for the number of samples, and a possible accuracy target if known upfront.\\
It would be worthwhile to further test the potential of MqPLHS in high-dimensional problems and with more engineering applications. Furthermore, FpPLHS and MqPLHS could be further improved by explicitly optimizing the projective properties to some extent. One focus of future research will certainly be aimed at comparing exploration-based and hybrid strategies for single and multi-response systems.
}

\section{Conflict of interest}\label{sec5}
On behalf of all authors, the corresponding author states that there is no conflict of interest.
\section{Funding Acknowledgments}\label{sec6}
The authors received no financial support for the research, authorship, and/or publication of this article.
\\
We would like to express our gratitude to Mr. Daniel Grealy for his valuable and constructive suggestions.

\bibliographystyle{unsrtnat} %%% unsrt, unsrtnat, apalike
\bibliography{references}  %%% Uncomment this line and comment out the ``thebibliography'' section below to use the external .bib file (using bibtex) .

\begin{thebibliography}{42}
\providecommand{\natexlab}[1]{#1}
\providecommand{\url}[1]{\texttt{#1}}
\expandafter\ifx\csname urlstyle\endcsname\relax
  \providecommand{\doi}[1]{doi: #1}\else
  \providecommand{\doi}{doi: \begingroup \urlstyle{rm}\Url}\fi

\bibitem[Baier et~al.(1994)Baier, See{\ss}elberg, and
  Specht]{Baier1994OptimierungID}
Horst Baier, Christoph See{\ss}elberg, and Bernhard Specht.
\newblock Optimierung in der strukturmechanik.
\newblock 1994.

\bibitem[Koch et~al.(2018)Koch, Mattern, and Bitsche]{Koch2008}
M.~Koch, S.~Mattern, and R.~Bitsche.
\newblock Facing future challenges in crash simulation engineering - model
  organization, quality and management at porsche.
\newblock \emph{15th International LS-DYNA Users Conference}, 2018.

\bibitem[Fukuda et~al.(2018)Fukuda, Pinto, dos Santos~Moreira, Saviano, and
  Lourenço]{Fukuda2018DesignOE}
Isa~Martins Fukuda, Camila Francini~Fidelis Pinto, Camila dos Santos~Moreira,
  Alessandro~Morais Saviano, and Felipe~Rebello Lourenço.
\newblock Design of experiments (doe) applied to pharmaceutical and analytical
  quality by design (qbd).
\newblock \emph{Brazilian Journal of Pharmaceutical Sciences}, 2018.

\bibitem[Kleijnen and van Beers(2004)]{Kleijnen2004}
J~P~C Kleijnen and W~C~M van Beers.
\newblock Application-driven sequential designs for simulation experiments:
  Kriging metamodelling.
\newblock \emph{Journal of the Operational Research Society}, 55\penalty0
  (8):\penalty0 876--883, 2004.
\newblock \doi{10.1057/palgrave.jors.2601747}.
\newblock URL \url{https://doi.org/10.1057/palgrave.jors.2601747}.

\bibitem[Johnson et~al.(1990)Johnson, Moore, and Ylvisaker]{JOHNSON1990131}
M.E. Johnson, L.M. Moore, and D.~Ylvisaker.
\newblock Minimax and maximin distance designs.
\newblock \emph{Journal of Statistical Planning and Inference}, 26\penalty0
  (2):\penalty0 131--148, 1990.
\newblock ISSN 0378-3758.
\newblock \doi{https://doi.org/10.1016/0378-3758(90)90122-B}.
\newblock URL
  \url{https://www.sciencedirect.com/science/article/pii/037837589090122B}.

\bibitem[Jiang et~al.(2015)Jiang, Shu, Zhou, Zhou, Shao, and Xu]{JIANG2015532}
Ping Jiang, Leshi Shu, Qi~Zhou, Hui Zhou, Xinyu Shao, and Junnan Xu.
\newblock A novel sequential exploration-exploitation sampling strategy for
  global metamodeling.
\newblock \emph{IFAC-PapersOnLine}, 48\penalty0 (28):\penalty0 532--537, 2015.
\newblock ISSN 2405-8963.
\newblock \doi{https://doi.org/10.1016/j .ifacol.2015.12.183}.
\newblock URL
  \url{https://www.sciencedirect.com/science/article/pii/S2405896315028074}.
\newblock 17th IFAC Symposium on System Identification SYSID 2015.

\bibitem[Crombecq et~al.(2011)Crombecq, Laermans, and Dhaene]{CROMBECQ2011683}
K.~Crombecq, E.~Laermans, and T.~Dhaene.
\newblock Efficient space-filling and non-collapsing sequential design
  strategies for simulation-based modeling.
\newblock \emph{European Journal of Operational Research}, 214\penalty0
  (3):\penalty0 683--696, 2011.
\newblock ISSN 0377-2217.
\newblock \doi{https://doi.org/10.1016/j .ejor.2011.05.032}.
\newblock URL
  \url{https://www.sciencedirect.com/science/article/pii/S0377221711004577}.

\bibitem[Lehmensiek et~al.(2002)Lehmensiek, Meyer, and Müller]{Lehmensiek2002}
Robert Lehmensiek, Petrie Meyer, and Martinette Müller.
\newblock Adaptive sampling applied to multivariate, multiple output rational
  interpolation models with application to microwave circuits.
\newblock \emph{International Journal of RF and Microwave Computer-Aided
  Engineering}, 12\penalty0 (4):\penalty0 332--340, 2002.
\newblock \doi{https://doi.org/10.1002/mmce.10032}.
\newblock URL \url{https://onlinelibrary.wiley.com/doi/abs/10.1002/mmce.10032}.

\bibitem[Sugiyama(2006)]{Sugiyama2006}
Masashi Sugiyama.
\newblock Active learning in approximately linear regression based on
  conditional expectation of generalization error.
\newblock \emph{J. Mach. Learn. Res.}, 7:\penalty0 141–166, dec 2006.
\newblock ISSN 1532-4435.

\bibitem[Crombecq et~al.(2009)Crombecq, Couckuyt, Gorissen, and
  Dhaene]{Crombecq2009SpacefillingSD}
Karel Crombecq, Ivo Couckuyt, Dirk Gorissen, and Tom Dhaene.
\newblock Space-filling sequential design strategies for adaptive surrogate
  modelling.
\newblock In \emph{SOCO 2009}, 2009.

\bibitem[Morris and Mitchell(1995)]{Morris1995ExploratoryDF}
Max~D. Morris and Toby~J. Mitchell.
\newblock Exploratory designs for computational experiments.
\newblock \emph{Journal of Statistical Planning and Inference}, 43:\penalty0
  381--402, 1995.

\bibitem[Zhang et~al.(2012)Zhang, Chowdhury, and Messac]{zhang2012adaptive}
Jie Zhang, Souma Chowdhury, and Achille Messac.
\newblock An adaptive hybrid surrogate model.
\newblock \emph{Structural and Multidisciplinary Optimization}, 46\penalty0
  (2):\penalty0 223--238, 2012.

\bibitem[Garud et~al.(2017)Garud, Karimi, and Kraft]{GARUD2017103}
Sushant~Suhas Garud, I.A. Karimi, and Markus Kraft.
\newblock Smart sampling algorithm for surrogate model development.
\newblock \emph{Computers \& Chemical Engineering}, 96:\penalty0 103--114,
  2017.
\newblock ISSN 0098-1354.
\newblock \doi{https://doi.org/10.1016/j .compchemeng.2016.10.006}.
\newblock URL
  \url{https://www.sciencedirect.com/science/article/pii/S0098135416303210}.

\bibitem[Husslage(2006)]{Husslage2006MaximinDF}
Bart Husslage.
\newblock Maximin designs for computer experiments.
\newblock 2006.

\bibitem[Draper and Pukelsheim(1996)]{Draper1996}
N.~R. Draper and F.~Pukelsheim.
\newblock An overview of design of experiments.
\newblock \emph{Statistical Papers}, 37:\penalty0 1--32, 1996.
\newblock \doi{10.1007/BF02926157}.
\newblock URL \url{https://doi.org/10.1007/BF02926157}.

\bibitem[Wu(2017)]{Wu2017}
C.~F.~Jeff Wu.
\newblock A fresh look at effect aliasing and interactions: some new wine in
  old bottles, 2017.
\newblock URL \url{https://arxiv.org/abs/1703.02113}.

\bibitem[van Dam et~al.(2007)van Dam, Husslage, den Hertog, and
  Melissen]{vanDam2007}
Edwin~R. van Dam, Bart Husslage, Dick den Hertog, and Hans Melissen.
\newblock Maximin latin hypercube designs in two dimensions.
\newblock \emph{Operations Research}, 55\penalty0 (1):\penalty0 158--169, 2007.
\newblock \doi{10.1287/opre.1060.0317}.
\newblock URL \url{https://doi.org/10.1287/opre.1060.0317}.

\bibitem[Sheikholeslami and Razavi(2017)]{SHEIKHOLESLAMI2017109}
Razi Sheikholeslami and Saman Razavi.
\newblock Progressive latin hypercube sampling: An efficient approach for
  robust sampling-based analysis of environmental models.
\newblock \emph{Environmental Modelling \& Software}, 93:\penalty0 109--126,
  2017.
\newblock ISSN 1364-8152.
\newblock \doi{https://doi.org/10.1016/j .envsoft.2017.03.010}.
\newblock URL
  \url{https://www.sciencedirect.com/science/article/pii/S1364815216305096}.

\bibitem[Grosso et~al.(2009)Grosso, Jamali, and Locatelli]{GROSSO2009541}
A.~Grosso, A.R.M.J.U. Jamali, and M.~Locatelli.
\newblock Finding maximin latin hypercube designs by iterated local search
  heuristics.
\newblock \emph{European Journal of Operational Research}, 197\penalty0
  (2):\penalty0 541--547, 2009.
\newblock ISSN 0377-2217.
\newblock \doi{https://doi.org/10.1016/j .ejor.2008.07.028}.
\newblock URL
  \url{https://www.sciencedirect.com/science/article/pii/S0377221708006644}.

\bibitem[Sobol'(1967)]{SOBOL196786}
I.M Sobol'.
\newblock On the distribution of points in a cube and the approximate
  evaluation of integrals.
\newblock \emph{USSR Computational Mathematics and Mathematical Physics},
  7\penalty0 (4):\penalty0 86--112, 1967.
\newblock ISSN 0041-5553.
\newblock \doi{https://doi.org/10.1016/0041-5553(67)90144-9}.
\newblock URL
  \url{https://www.sciencedirect.com/science/article/pii/0041555367901449}.

\bibitem[Halton and Smith(1964)]{Halton1964Algorithm2R}
John~H. Halton and G.~B. Smith.
\newblock Algorithm 247: Radical-inverse quasi-random point sequence.
\newblock \emph{Commun. ACM}, 7:\penalty0 701--702, 1964.

\bibitem[Hammersley and Handscomb(1964)]{Hammersley1964}
J.~M. Hammersley and D.~C. Handscomb.
\newblock \emph{The General Nature of Monte Carlo Methods}, pages 1--9.
\newblock Springer Netherlands, Dordrecht, 1964.
\newblock ISBN 978-94-009-5819-7.
\newblock \doi{10.1007/978-94-009-5819-7_1}.
\newblock URL \url{https://doi.org/10.1007/978-94-009-5819-7_1}.

\bibitem[Saltelli et~al.(2010)Saltelli, Annoni, Azzini, Campolongo, Ratto, and
  Tarantola]{SALTELLI2010259}
Andrea Saltelli, Paola Annoni, Ivano Azzini, Francesca Campolongo, Marco Ratto,
  and Stefano Tarantola.
\newblock Variance based sensitivity analysis of model output. design and
  estimator for the total sensitivity index.
\newblock \emph{Computer Physics Communications}, 181\penalty0 (2):\penalty0
  259--270, 2010.
\newblock ISSN 0010-4655.
\newblock \doi{https://doi.org/10.1016/j .cpc.2009.09.018}.
\newblock URL
  \url{https://www.sciencedirect.com/science/article/pii/S0010465509003087}.

\bibitem[Burhenne et~al.(2011)Burhenne, Jacob, and
  Henze]{Burhenne2011SAMPLINGBO}
Sebastian Burhenne, Dirk Jacob, and Gregor~P. Henze.
\newblock Sampling based on sobol 0 sequences for monte carlo techniques
  applied to building simulations.
\newblock 2011.

\bibitem[Du et~al.(1999)Du, Faber, and Gunzburger]{Du1999}
Qiang Du, Vance Faber, and Max Gunzburger.
\newblock Centroidal voronoi tessellations: Applications and algorithms.
\newblock \emph{SIAM Review}, 41\penalty0 (4):\penalty0 637--676, 1999.
\newblock \doi{10.1137/S0036144599352836}.
\newblock URL \url{https://doi.org/10.1137/S0036144599352836}.

\bibitem[Burns(2009)]{Burns2009CENTROIDALVT}
Jared Burns.
\newblock Centroidal voronoi tessellations.
\newblock 2009.

\bibitem[Feinberg and Langtangen(2015)]{FEINBERG201546}
Jonathan Feinberg and Hans~Petter Langtangen.
\newblock Chaospy: An open source tool for designing methods of uncertainty
  quantification.
\newblock \emph{Journal of Computational Science}, 11:\penalty0 46--57, 2015.
\newblock ISSN 1877-7503.
\newblock \doi{https://doi.org/10.1016/j .jocs.2015.08.008}.
\newblock URL
  \url{https://www.sciencedirect.com/science/article/pii/S1877750315300119}.

\bibitem[Schreiber et~al.(2006)Schreiber, Nora, Stage, Barlow, and
  King]{Schreiber2006}
James~B. Schreiber, Amaury Nora, Frances~K. Stage, Elizabeth~A. Barlow, and
  Jamie King.
\newblock Reporting structural equation modeling and confirmatory factor
  analysis results: A review.
\newblock \emph{The Journal of Educational Research}, 99\penalty0 (6):\penalty0
  323--338, 2006.
\newblock \doi{10.3200/JOER.99.6.323-338}.
\newblock URL \url{https://doi.org/10.3200/JOER.99.6.323-338}.

\bibitem[Jones et~al.(1998)Jones, Schonlau, and Welch]{Jones1998EfficientGO}
Donald~R. Jones, Matthias Schonlau, and William~J. Welch.
\newblock Efficient global optimization of expensive black-box functions.
\newblock \emph{Journal of Global Optimization}, 13:\penalty0 455--492, 1998.

\bibitem[Winston(1992)]{Winston1992}
Patrick~Henry Winston.
\newblock \emph{Artificial Intelligence (3rd Ed.)}.
\newblock Addison-Wesley Longman Publishing Co., Inc., USA, 1992.
\newblock ISBN 0201533774.

\bibitem[Salomon(1996)]{SALOMON1996263}
Ralf Salomon.
\newblock Re-evaluating genetic algorithm performance under coordinate rotation
  of benchmark functions. a survey of some theoretical and practical aspects of
  genetic algorithms.
\newblock \emph{Biosystems}, 39\penalty0 (3):\penalty0 263--278, 1996.
\newblock ISSN 0303-2647.
\newblock \doi{https://doi.org/10.1016/0303-2647(96)01621-8}.
\newblock URL
  \url{https://www.sciencedirect.com/science/article/pii/0303264796016218}.

\bibitem[Garc{\'{\i}}a{-}Pedrajas et~al.(2011)Garc{\'{\i}}a{-}Pedrajas,
  Herv{\'{a}}s{-}Mart{\'{\i}}nez, and Ortiz{-}Boyer]{Boyer2011}
Nicol{\'{a}}s Garc{\'{\i}}a{-}Pedrajas, C{\'{e}}sar
  Herv{\'{a}}s{-}Mart{\'{\i}}nez, and Domingo Ortiz{-}Boyer.
\newblock {CIXL2:} {A} crossover operator for evolutionary algorithms based on
  population features.
\newblock \emph{CoRR}, abs/1109.2146, 2011.
\newblock URL \url{http://arxiv.org/abs/1109.2146}.

\bibitem[Jamil and Yang(2013)]{Jamil2013ALS}
Momin Jamil and Xin-She Yang.
\newblock A literature survey of benchmark functions for global optimisation
  problems.
\newblock \emph{Int. J. Math. Model. Numer. Optimisation}, 4:\penalty0
  150--194, 2013.

\bibitem[Adorio and January(2005)]{Adorio2005MVFM}
Ernesto~Padernal Adorio and Revised January.
\newblock Mvf - multivariate test functions library in c for unconstrained
  global optimization.
\newblock 2005.

\bibitem[Molga and Smutnicki(2005)]{Molga2005}
M.~Molga and C.~Smutnicki.
\newblock Test functions for optimization needs.
\newblock 2005.
\newblock URL \url{http://www.zsd.ict.pwr.wroc.pl/files/docs/functions.pdf}.

\bibitem[Yao et~al.(2003)Yao, Liu, Liang, and Lin]{Yao2003}
Xin Yao, Yong Liu, Ko-Hsin Liang, and Guangming Lin.
\newblock \emph{Fast Evolutionary Algorithms}, pages 45--94.
\newblock Springer Berlin Heidelberg, Berlin, Heidelberg, 2003.
\newblock ISBN 978-3-642-18965-4.
\newblock \doi{10.1007/978-3-642-18965-4_2}.
\newblock URL \url{https://doi.org/10.1007/978-3-642-18965-4_2}.

\bibitem[Duvenaud(2014)]{duvenaud-thesis-2014}
David Duvenaud.
\newblock \emph{Automatic Model Construction with {G}aussian Processes}.
\newblock PhD thesis, {Computational and Biological Learning Laboratory,
  University of Cambridge}, 2014.

\bibitem[Rasmussen and Williams(2006)]{Rasmussen2006}
Carl~Eduard Rasmussen and Christopher~K.I. Williams.
\newblock Gaussian processes for machine learning.
\newblock 2006.

\bibitem[Brereton and Lloyd(2010)]{Brereton2010SupportVM}
Richard~G. Brereton and Gavin~Rhys Lloyd.
\newblock Support vector machines for classification and regression.
\newblock \emph{The Analyst}, 135 2:\penalty0 230--67, 2010.

\bibitem[Gunn(1998)]{Gunn1998SupportVM}
Steve~R. Gunn.
\newblock Support vector machines for classification and regression.
\newblock 1998.

\bibitem[Verleysen and François(2005)]{Verleysen2005TheCO}
Michel Verleysen and Damien François.
\newblock The curse of dimensionality in data mining and time series
  prediction.
\newblock In \emph{IWANN}, 2005.

\bibitem[Forrester et~al.(2008)Forrester, Sóbester, and Keane]{Forrester2008}
A.I.J. Forrester, A.~Sóbester, and A.J. Keane.
\newblock \emph{Sampling Plans}, chapter~1, pages 1--31.
\newblock John Wiley \& Sons, Ltd, 2008.
\newblock \doi{https://doi.org/10.1002/9780470770801.ch1}.
\newblock URL
  \url{https://onlinelibrary.wiley.com/doi/abs/10.1002/9780470770801.ch1}.

\end{thebibliography}

%%% Uncomment this section and comment out the \bibliography{references} line above to use inline references.
% \begin{thebibliography}{1}

% 	\bibitem{kour2014real}
% 	George Kour and Raid Saabne.
% 	\newblock Real-time segmentation of on-line handwritten arabic script.
% 	\newblock In {\em Frontiers in Handwriting Recognition (ICFHR), 2014 14th
% 			International Conference on}, pages 417--422. IEEE, 2014.

% 	\bibitem{kour2014fast}
% 	George Kour and Raid Saabne.
% 	\newblock Fast classification of handwritten on-line arabic characters.
% 	\newblock In {\em Soft Computing and Pattern Recognition (SoCPaR), 2014 6th
% 			International Conference of}, pages 312--318. IEEE, 2014.

% 	\bibitem{hadash2018estimate}
% 	Guy Hadash, Einat Kermany, Boaz Carmeli, Ofer Lavi, George Kour, and Alon
% 	Jacovi.
% 	\newblock Estimate and replace: A novel approach to integrating deep neural
% 	networks with existing applications.
% 	\newblock {\em arXiv preprint arXiv:1804.09028}                                   , 2018.

% \end{thebibliography}

\end{document}